\documentclass{article}

\usepackage{ifthen}

\newcommand{\arxiv}{}

\usepackage{amsmath,amsfonts,amssymb,algorithm,algorithmic,psfrag,amsthm}
\usepackage{graphicx}
\usepackage{array}
\usepackage{psfrag}
\usepackage{subfigure}
\usepackage{bm}
\usepackage{color}
\usepackage[curve]{xypic}
\usepackage{multirow}
\usepackage[hidelinks]{hyperref}

\ifthenelse{\isundefined{\arxiv}}{}{\usepackage{fullpage}}

\newtheorem{theorem}{Theorem}[section]
\newtheorem{lemma}[theorem]{Lemma}
\newtheorem{corollary}[theorem]{Corollary}

\theoremstyle{definition}
\newtheorem{definition}[theorem]{Definition}

\theoremstyle{remark}
\newtheorem{remark}[theorem]{Remark}
\numberwithin{equation}{section}

\newcommand{\A} {{\mathbb A}}
\newcommand{\B} {{\mathbb B}}

\newcommand{\R} {\mathbb R}
\newcommand{\Ell} {\mathbb L}

\newcommand{\p}{\partial}

\newcommand{\cancel}[1]{}
\newcommand{\raw}{\rightarrow}

\newcommand{\eps}{\epsilon}

\newcommand{\vn}[1]{\left|\left|#1\right|\right|}
\newcommand{\vsn}[1]{\left|#1\right|}

\newcommand{\ds}{\displaystyle}

\newcommand{\ltri}{\lambda^{\rm Tri}}
\newcommand{\lwach}{\lambda^{\rm Wach}}
\newcommand{\lopt}{\lambda^{\rm Har}}
\newcommand{\lsib}{\lambda^{\rm Sibs}}
\newcommand{\lmval}{\lambda^{\rm MVal}}

\newcommand\diam{\textnormal{diam}}

\newcommand{\bc}{\textbf{c}}
\newcommand{\bp}{\textbf{p}}

\newcommand{\bv}{\textbf{v}}
\newcommand{\bx}{\textbf{x}}

\newcommand{\bal}{{\bf \alpha}}

\newcommand{\lpn}[3]{\vn{#1}_{L^{#2}(#3)}}
\newcommand{\hpn}[3]{\vn{#1}_{H^{#2}(#3)}}   
\newcommand{\hpsn}[3]{\vsn{#1}_{H^{#2}(#3)}} 

\newcommand{\cP}{{\mathcal P}}
\newcommand{\cS}{{\mathcal S}}

\newcommand{\id}{\mathbb I}
\newcommand{\Qx}{\text{Q}_{\xi}}
\newcommand{\Qc}{\text{Q}_{c}}
\newcommand{\tps}{(\theta+\sigma)}
\newcommand{\tms}{(\theta-\sigma)}

\newtheorem{proposition}[theorem]{Proposition}

\newcommand{\comment}[1]{}

\begin{document}

\ifthenelse{\isundefined{\arxiv}}{
\title[Quadratic Serendipity Finite Elements on Polygons]{Quadratic Serendipity Finite Elements on Polygons Using Generalized Barycentric Coordinates}
}{
\title{Quadratic Serendipity Finite Elements on Polygons Using Generalized Barycentric Coordinates\thanks{This research was supported in part by NIH contracts R01-EB00487, R01-GM074258, and a grant from the UT-Portugal CoLab project.}}
}


\ifthenelse{\isundefined{\arxiv}}{
\author{Alexander Rand}
\address{Institute for Computational Engineering and Sciences, The University of Texas at Austin}
\email{arand@ices.utexas.edu}

\author{Andrew Gillette}
\address{Department of Mathematics, The University of Texas at Austin}
\email{agillette@math.utexas.edu}

\author{Chandrajit Bajaj}
\address{Department of Computer Science, The University of Texas at Austin}
\email{bajaj@cs.utexas.edu}

\thanks{This research was supported in part by NIH contracts R01-EB00487, R01-GM074258, and a
grant from the UT-Portugal CoLab project.}
}{

\author{Alexander Rand\thanks{Institute for Computational Engineering and Sciences, The University of Texas at Austin, \href{mailto:arand@ices.utexas.edu}{arand@ices.utexas.edu}}, Andrew Gillette\thanks{Department of Mathematics, The University of Texas at Austin, \href{mailto:agillette@math.utexas.edu}{agillette@math.utexas.edu}}, and Chandrajit Bajaj\thanks{Department of Computer Science, The University of Texas at Austin, \href{mailto:bajaj@cs.utexas.edu}{bajaj@cs.utexas.edu}}}

}


\ifthenelse{\isundefined{\arxiv}}{

\date{}

\dedicatory{}

}{}

\maketitle

\begin{abstract}
We introduce a finite element construction for use on the class of convex, planar polygons and show it obtains a quadratic error convergence estimate.  
On a convex $n$-gon satisfying simple geometric criteria, our construction produces $2n$ basis functions, associated in a Lagrange-like fashion to each vertex and each edge midpoint, by transforming and combining a set of $n(n+1)/2$ basis functions known to obtain quadratic convergence.  
The technique broadens the scope of the so-called `serendipity' elements, previously studied only for quadrilateral and regular hexahedral meshes, by employing the theory of generalized barycentric coordinates.
Uniform \textit{a priori} error estimates are established over the class of convex quadrilaterals with bounded aspect ratio as well as over the class of generic convex planar polygons satisfying additional shape regularity conditions to exclude large interior angles and short edges. 
Numerical evidence is provided on a trapezoidal quadrilateral mesh, previously not amenable to serendipity constructions, and applications to adaptive meshing are discussed. 
\end{abstract}

\section{Introduction}

Barycentric coordinates provide a basis for linear finite elements on simplices, and generalized barycentric coordinates naturally produce a suitable basis for linear finite elements on general polygons.  Various applications make use of this technique~\cite{GB2010,GB2011,MKBWG2008,MP2008,RS2006,SAB2010,SM2006,ST2004,TS06,WBG07}, but in each case, only linear error estimates can be asserted.  A quadratic finite element can easily be constructed by taking pairwise products of the basis functions from the linear element, yet this approach has not been pursued, primarily since the requisite number of basis functions grows quadratically in the number of vertices of the polygon.  Still, many of the pairwise products are zero along the entire polygonal boundary and thus are unimportant for inter-element continuity, a key ingredient in finite element theory.  For quadrilateral elements, these `extra' basis functions are well understood and, for quadrilaterals that can be affinely mapped to a square, the so-called `serendipity element' yields an acceptable basis consisting of only those basis functions needed to guarantee inter-element continuity~\cite{ZT2000,ABF02,AA10}.  We generalize this construction to produce a quadratic serendipity element for arbitrary convex polygons derived from generalized barycentric coordinates.

Our construction yields a set of Lagrange-like basis functions $\{\psi_{ij}\}$ -- one per vertex and one per edge midpoint -- using a linear combination of pairwise products of generalized barycentric functions $\{\lambda_i\}$.  
We show that this set spans all constant, linear, and quadratic polynomials, making it suitable for finite element analysis via the Bramble-Hilbert lemma.  
Further, given uniform bounds on the aspect ratio, edge length, and interior angles of the polygon, we bound $\hpn{\psi_{ij}}{1}{\Omega}$ uniformly with respect to $\hpn{\lambda_i}{1}{\Omega}$.
Since our previous work shows that $\hpn{\lambda_i}{1}{\Omega}$ is bounded uniformly under these geometric hypotheses, this proves that the $\psi_{ij}$ functions are well-behaved.

\begin{figure}[t]
\xymatrix @R=.02in{
\parbox{.18\textwidth}{\includegraphics[width=.18\textwidth]{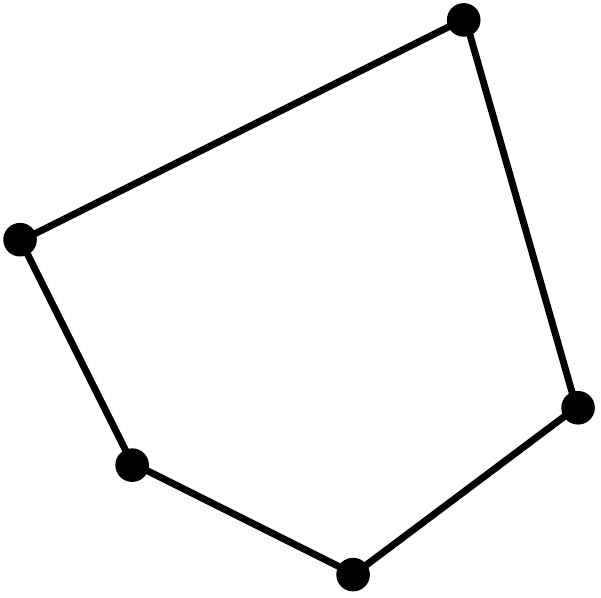}} &
\parbox{.18\textwidth}{\includegraphics[width=.18\textwidth]{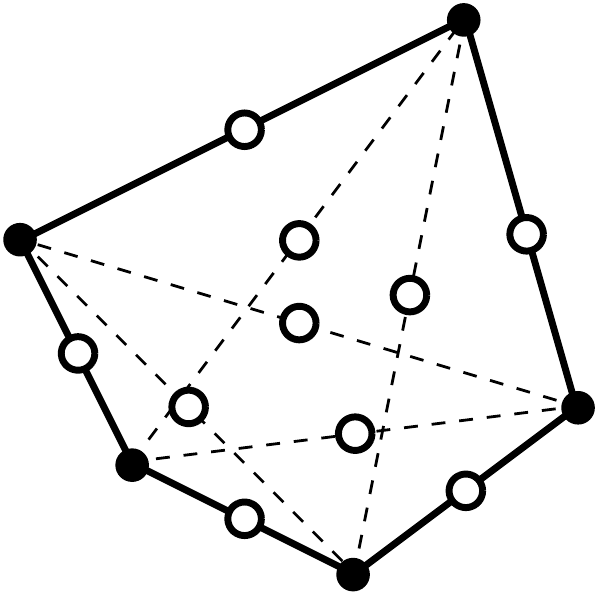}} &
\parbox{.18\textwidth}{\includegraphics[width=.18\textwidth]{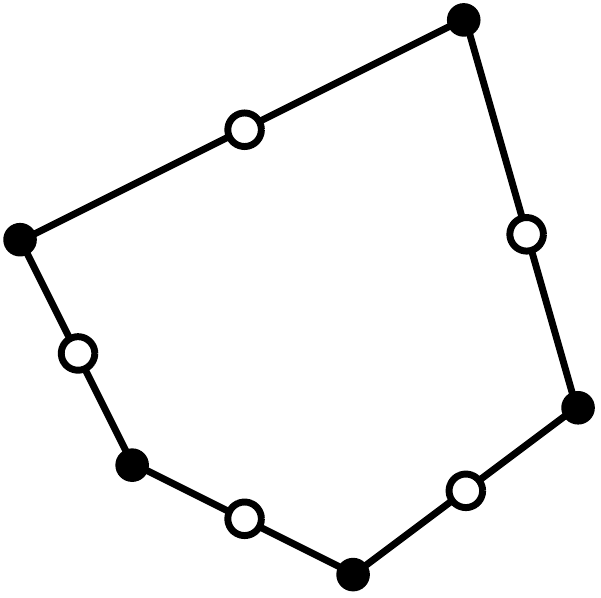}} &
\parbox{.18\textwidth}{\includegraphics[width=.18\textwidth]{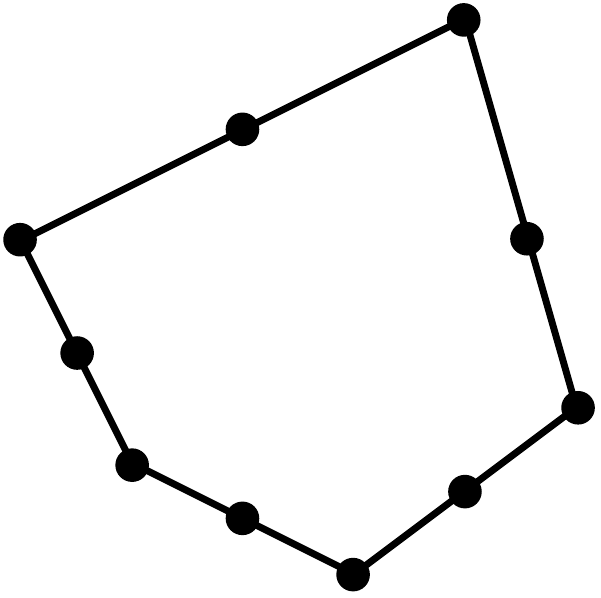}} \\
{\{\lambda_i\}} \ar[r]^-{\text{pairwise}}_-{\text{products}} &  {\{\mu_{ab}\}}   \ar[r]^{\A} & {\{\xi_{ij}\}}  \ar[r]^{\B} & {\{\psi_{ij}\}} \\
\text{Linear} & \text{Quadratic} & \text{Serendipity} & \text{Lagrange}
}
\caption{Overview of the construction process.  In each figure, the dots are in one-to-one correspondence with the set of functions listed below it.  At filled dots, all functions in the set evaluate to zero except for the function corresponding to the dot which evaluates to one.  The rightmost element has quadratic precision with only these types of `Lagrange-like' basis functions.
}
\label{fg:motiv}
\end{figure}

Figure~\ref{fg:motiv} gives a visual depiction of the construction process.  Starting with one generalized barycentric function $\lambda_i$ per vertex of an $n$-gon, take all pairwise products yielding a total of $n(n+1)/2$ functions $\mu_{ab}:=\lambda_a\lambda_b$.  The linear transformation $\A$ reduces the set $\{\mu_{ab}\}$ to the $2n$ element set $\{\xi_{ij}\}$, indexed over vertices and edge midpoints of the polygon.  A simple bounded linear transformation $\B$ converts $\{\xi_{ij}\}$ into a basis $\{\psi_{ij}\}$ which satisfies the ``Lagrange property'' meaning each function takes the value 1 at its associated node and 0 at all other nodes.

The paper is organized as follows.  
In Section~\ref{sec:bkgdnot} we review relevant background on finite element theory, serendipity elements, and generalized barycentric functions.  
In Section~\ref{sec:construct}, we show that if the entries of matrix $\A$ satisfy certain linear constraints $\Qc1$-$\Qc3$, the resulting set of functions $\{\xi_{ij}\}$ span all constant, linear and quadratic monomials in two variables, a requirement for quadratic finite elements.
In Section~\ref{sec:specialcase}, we show how the constraints $\Qc1$-$\Qc3$ can be satisfied in the special cases of the unit square, regular polygons, and convex quadrilaterals.  
In Section~\ref{sec:genpoly}, we show how $\Qc1$-$\Qc3$ can be satisfied on a simple convex polygon.
We also prove that the resulting value of $||\A||$ is bounded uniformly, provided the convex polygon satisfies certain geometric quality conditions.
In Section~\ref{sec:lag-reduc} we define $\B$ and show that the final $\{\psi_{ij}\}$ basis is Lagrange-like.  
Finally, in Section~\ref{sec:conc}, we describe practical applications, give numerical evidence, and consider future directions.

\section{Background and Notation}
\label{sec:bkgdnot}

\begin{figure}
\centering
\includegraphics[width=.35\textwidth]{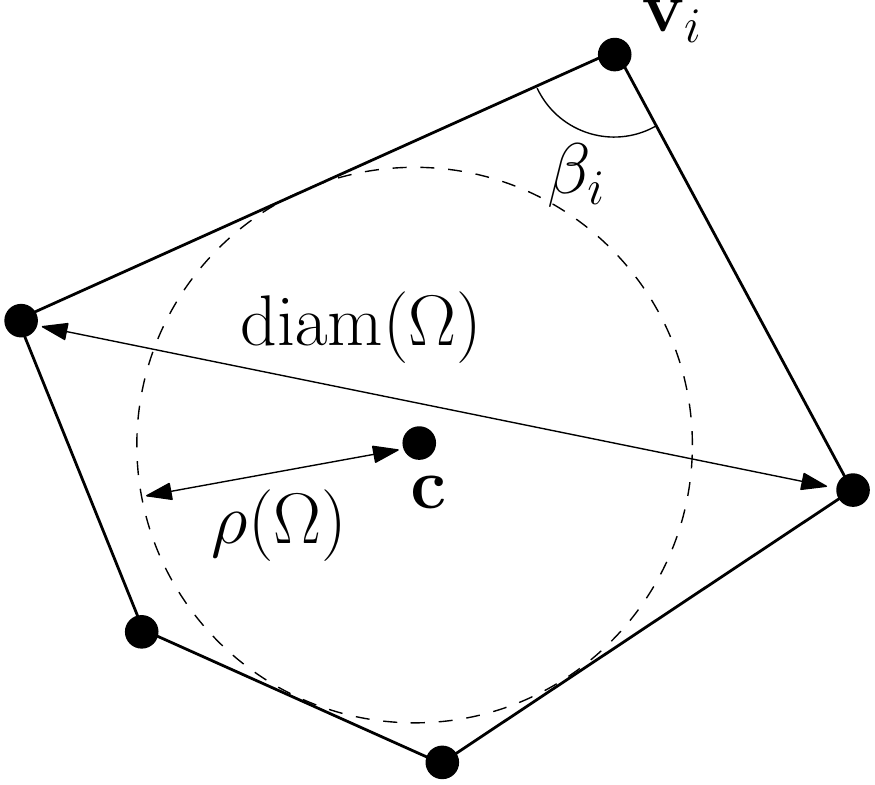}
\caption{Notation used to describe polygonal geometry.}\label{fg:aspectRatio}
\end{figure}

Let $\Omega$ be a convex polygon with $n$ vertices $(\bv_1, \ldots, \bv_n)$ ordered counter-clockwise.
Denote the interior angle at $\bv_i$ by $\beta_i$.
The largest distance between two points in $\Omega$ (the diameter of $\Omega$) is denoted $\diam(\Omega)$ and the radius of the largest inscribed circle is denoted $\rho(\Omega)$.
The center of this circle is denoted $\bc$ and is selected arbitrarily when no unique circle exists.
The \textbf{aspect ratio} (or chunkiness parameter) $\gamma$ is the ratio of the diameter to the radius of the largest inscribed circle, i.e.
\[\gamma := \frac{\diam(\Omega)}{\rho(\Omega)}.\]
The notation is shown in Figure~\ref{fg:aspectRatio}.


For a multi-index $\bal = (\alpha_1, \alpha_2)$ and point $\bx = (x,y)$, define $\bx^\bal := x^{\alpha_1} y^{\alpha_2}$, $\alpha ! := \alpha_1 \alpha_2$, $|\bal| := \alpha_1 + \alpha_2$, and $D^\bal u := \p^{|\bal|} u/\p x^{\alpha_1}\p y^{\alpha_2}$.
The Sobolev semi-norms and norms over an open set $\Omega$ are defined by
\begin{align*}
\hpsn{u}{m}{\Omega}^2 &:=  \int_\Omega \sum_{|\alpha| = m} |D^\alpha u(\bx)|^2 \,{\rm d} \bx &{\rm and} & & \hpn{u}{m}{\Omega}^2 &:= \sum_{0\leq k\leq m}\hpsn{u}{k}{\Omega}^2.
\end{align*}
The $H^{0}$-norm is the $L^2$-norm and will be denoted $\lpn{\cdot}{2}{\Omega}$.  The space of polynomials of degree $\leq k$ on a domain is denoted $\cP_k$.

\subsection{The Bramble-Hilbert Lemma}

A finite element method approximates a function $u$ from an infinite-dimensional functional space $V$ by a function $u_h$ from a finite-dimensional subspace $V_h\subset V$.  
One goal of such approaches is to prove that the error of the numerical solution $u_h$ is bounded \textit{a priori} by the error of the best approximation available in $V_h$, i.e. $\vn{u-u_h}_{V} \leq C \inf_{w\in V_h} \vn{u-w}_{V}$.
In this paper, $V=H^1$ and $V_h$ is the span of a set of functions defined piecewise over a 2D mesh of convex polygons.  The parameter $h$ indicates the maximum diameter of an element in the mesh.  Further details on the finite element method can be found in a number of textbooks~\cite{Ci02,BS08,EG04,ZT2000}.

A quadratic finite element method in this context means that when $h\raw 0$, the best approximation error ($\inf_{w\in V_h} \vn{u-w}_{V}$) converges to zero with order $h^2$.  This means the space $V_h$ is `dense enough' in $V$ to allow for quadratic convergence.  Such arguments are usually proved via the Bramble-Hilbert lemma which guarantees that if $V_h$ contains polynomials up to a certain degree, a bound on the approximation error can be found. The variant of the Bramble-Hilbert lemma stated below includes a uniform constant over all convex domains which is a necessary detail in the context of general polygonal elements and generalized barycentric functions. 

\begin{lemma}[Bramble-Hilbert~\cite{Ve99,DL04}]\label{lem:bramblehilbert}
There exists a uniform constant $C_{BH}$ such that for all convex polygons $\Omega$ and  for all $u\in H^{k}(\Omega)$, there exists a degree $k$ polynomial $p_u$ with
$\hpn{u-p_u}{k'}{\Omega} \leq C_{BH} \,\diam(\Omega)^{k+1-k'} \hpsn{u}{k+1}{\Omega}$ for any $k'\leq k$.
\end{lemma}

Our focus is on quadratic elements (i.e., $k=2$) and error estimates in the $H^1$-norm (i.e., $k'=1$) which yields an estimate that scales with $\diam(\Omega)^{2}$.  
Our methods extend to more general Sobolev spaces (i.e., $W^{k,p}$, the space of functions with all derivatives of order $\leq k$ in $L^{p}$) whenever the Bramble-Hilbert lemma holds.
Extensions to higher order elements ($k>2$) will be briefly discussed in Section~\ref{sec:conc}.

\begin{figure}
\centering
\includegraphics[width=.5\textwidth]{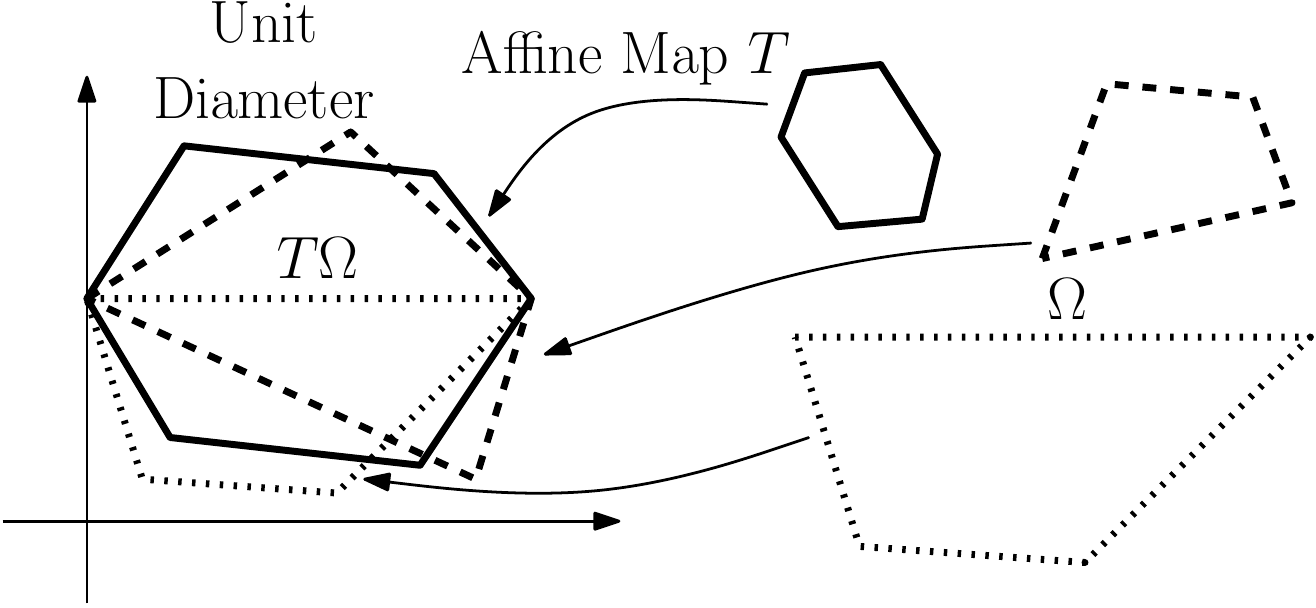}
\caption{Using affine transformation, analysis can be restricted to a class of unit diameter polygons.}\label{fg:affinetrans}
\end{figure}

Observe that if $\Omega$ is transformed by any invertible affine map $T$, the polynomial $p\circ T^{-1}$ on $T\Omega$ has the same degree as the polynomial $p$ on $\Omega$.  
This fact is often exploited in the simpler and well-studied case of triangular meshes; an estimate on a reference triangle $\hat K$ becomes an estimate on any physical triangle $K$ by passing through an affine transformation taking $\hat K$ to $K$.  
For $n>3$, however, two generic $n$-gons may differ by a non-affine transformation and thus, as we will see in the next section, the use of a single reference element can become overly restrictive on element geometry.  
In our arguments, we instead analyze classes of ``reference'' elements, namely, diameter one convex quadrilaterals or convex polygons of diameter one satisfying the geometric criteria given in Section~\ref{subsec:genbary}; see Figure~\ref{fg:affinetrans}.
Using a class of reference elements allows us to establish uniform error estimates over all affine transformations of this class.

\subsection{Serendipity Quadratic Elements}

The term `serendipity element' refers to a long-standing observation in the finite element community that tensor product bases of polynomials on rectangular meshes of quadrilaterals in 2D or cubes in 3D can obtain higher order convergence rates with fewer than the `expected' number of basis functions resulting from tensor products.  This phenomenon is discussed in many finite element textbooks, e.g.~\cite{SF73,H1987,Ci02}, and was recently characterized precisely by Arnold and Awanou~\cite{AA10}.  For instance, the degree $r$ tensor product basis on a square reference element has $(r+1)^2$ basis functions and can have guaranteed convergence rates of order $r+1$ when transformed to a rectangular mesh via bilinear isomorphisms~\cite{ABF02}.  By the Bramble-Hilbert lemma, however, the function space spanned by this basis may be unnecessarily large as the dimension of $\cP_r$ is only $(r+1)(r+2)/2$ and only $4r$ degrees of freedom associated to the boundary are needed to ensure sufficient inter-element continuity in $H^1$.

This motivates the construction of the serendipity element for quadrilaterals.  
By a judicious choice of basis functions, an order $r$ convergence rate can be obtained with one basis function associated to each vertex, $(r-1)$ basis functions associated to each edge, and $q$ additional functions associated to interior points of the quadrilateral, where $q=0$ for $r<4$ and $q=(r-2)(r-1)/2$ for $r\geq 4$~\cite{AA10}.  
Such an approach only works if the reference element is mapped via an affine transformation; it has been demonstrated that the serendipity element fails on trapezoidal elements, such as those shown in Figure~\ref{fg:nonaffinequad}~\cite{MH92,KOF99,ZT2000,ZK00}.

Some very specific serendipity elements have been constructed for quadrilaterals and regular hexagons based on the Wachspress coordinates (discussed in the next sections)~\cite{W1975,apprato1979elements,gout1979construction,apprato1979rational,gout1985rational}.
Our work generalizes this construction to arbitrary polygons without dependence on the type of generalized barycentric coordinate selected and  with uniform bounds under certain geometric criteria.

\subsection{Generalized Barycentric Elements}
\label{subsec:genbary}

To avoid non-affine transformations associated with tensor products constructions on a single reference element, we use generalized barycentric coordinates to define our basis functions.
These coordinates are any functions satisfying the following agreed-upon definition in the literature.


\begin{definition}\label{def:barcoor}
Functions $\lambda_i:\Omega\raw\R$, $i=1,\ldots, n$ are \textbf{barycentric coordinates} on $\Omega$ if they satisfy two properties.
\renewcommand{\labelenumi}{B\arabic{enumi}.}
\begin{enumerate}
\item \textbf{Non-negative}:  $\lambda_i\geq 0$ on $\Omega$.\label{b:nonneg}
\item \textbf{Linear Completeness}: For any linear function $L:\Omega\raw\R$,
$\ds L=\sum_{i=1}^{n} L(\bv_i)\lambda_i$.\label{b:lincomp}
\end{enumerate}
\end{definition}


We will further restrict our attention to barycentric coordinates satisfying the following invariance property.
Let $T:\R^2 \rightarrow \R^2$ be a composition of translation, rotation, and uniform scaling transformations and let $\{\lambda^T_i\}$ denote a set of barycentric coordinates on $T\Omega$.

\renewcommand{\labelenumi}{B\arabic{enumi}.}
\begin{enumerate}
\setcounter{enumi}{2}
\item \textbf{Invariance:} $\ds\lambda_i(\bx)=\lambda_i^T(T(\bx))$.\label{b:invariance}
\end{enumerate}

This assumption will allow estimates over the class of convex sets with diameter one to be immediately extended to generic sizes since translation, rotation and uniform scaling operations can be easily passed through Sobolev norms.  At the expense of requiring uniform bounds over a class of diameter-one domains rather than a single reference element, we avoid having to handle non-affine mappings between reference and physical elements.

A set of barycentric coordinates $\{\lambda_i\}$ also satisfies three additional familiar properties.  A proof that B\ref{b:nonneg} and B\ref{b:lincomp} imply the additional properties B\ref{b:partition}-B\ref{b:interpolation} can be found in~\cite{GRB2011}.
Note that B\ref{b:partition} and B\ref{b:linprec} follow immediately by setting $L=1$ or $L=\bx$ in B\ref{b:lincomp}.
\renewcommand{\labelenumi}{B\arabic{enumi}.}
\begin{enumerate}
\setcounter{enumi}{3}
\item \textbf{Partition of unity:} $\ds\sum_{i=1}^{n}\lambda_i\equiv 1$. \label{b:partition}
\item \textbf{Linear precision:} $\ds\sum_{i=1}^{n}\bv_i\lambda_i(\bx)=\bx$. \label{b:linprec}
\item \textbf{Interpolation:} $\ds\lambda_i(\bv_j) = \delta_{ij}$. \label{b:interpolation}
\end{enumerate}

Various particular barycentric coordinates have been constructed in the literature.  We briefly mention a few of the more prominent kinds and associated references here; readers are referred to our prior work~\cite[Section 2]{GRB2011} as well as the survey papers of Cueto et al. \cite{CSCMCD2003} and Sukumar and Tabarraei \cite{ST2004} for further details.
The triangulation coordinates $\ltri$ are defined by triangulating the polygon and using the standard barycentric coordinates over each triangle~\cite{FHK2006}.  Harmonic coordinates $\lopt$ are defined as the solution to Laplace's equation on the polygon with piecewise linear boundary data satisfying B\ref{b:interpolation}~\cite{JMRGS07,MKBWG2008,C2008}.  Explicitly constructed functions include the rational Wachspress coordinates $\lwach$~\cite{W1975}, the Sibson coordinates $\lsib$ defined in terms of the Voronoi diagram of the vertices of the polygon~\cite{S1980,F1990}, and the mean value coordinates $\lmval$ defined by Floater~\cite{F2003,FHK2006}.

To obtain convergence estimates with any of these functions, certain geometric conditions must be satisfied by a generic mesh element.  We will consider domains satisfying the following three geometric conditions.

\renewcommand{\labelenumi}{G\arabic{enumi}.}
\begin{enumerate}
\item \textbf{Bounded aspect ratio:} There exists $\gamma^*\in\R$ such that $\gamma < \gamma^*$.  \label{g:ratio}
\item \textbf{Minimum edge length: } There exists $d_*\in\R$ such that $|\bv_i - \bv_j| > d_* > 0$ for all $i\neq j$. \label{g:minedge}
\item \textbf{Maximum interior angle:} There exists $\beta^*\in\R$ such that $\beta_i < \beta^* < \pi$ for all $i$.\label{g:maxangle}
\end{enumerate}

Under some set of these conditions, the $H^1$-norm of many generalized barycentric coordinates are bounded in $H^1$ norm.
This is a key estimate in asserting the expected (linear) convergence rate in the typical finite element setting.

\begin{theorem}[\cite{RGB2011b} for $\lmval$ and \cite{GRB2011} for others]\label{th:aicmresult}
For any convex polygon $\Omega$ satisfying G1, G2, and G3, $\ltri$, $\lopt$, $\lwach$, $\lsib$, and $\lmval$
are all bounded in $H^1$, i.e. there exists a constant $C>0$ such that
\begin{equation}\label{eq:aicmresult}
\hpn{\lambda_i}{1}{\Omega} \leq C.
\end{equation}
\end{theorem}

The results in \cite{GRB2011} and \cite{RGB2011b} are somewhat stronger than the statement of Theorem~\ref{th:aicmresult}, namely, not all of the geometric hypotheses are necessary for every coordinate type.
Our results, however, rely generically on any set of barycentric coordinates satisfying (\ref{eq:aicmresult}).
Any additional dependence on the shape geometry will be made explicitly clear in the proofs.
Weakening of the geometric hypotheses is discussed in Section~\ref{sec:conc}.

\subsection{Quadratic Precision Barycentric Functions}

Since generalized barycentric coordinates are only guaranteed to have linear precision (property B\ref{b:linprec}), they cannot provide greater than linear order error estimates.
Pairwise products of barycentric coordinates, however, provide quadratic precision as the following simple proposition explains.

\begin{proposition}
\label{prop:pair-prod}
Given a set of barycentric coordinates $\{\lambda_i\}_{i=1}^n$, the set of functions $\{\mu_{ab}\} := \{\lambda_a\lambda_b\}_{a,b=1}^n$ has constant, linear, and quadratic precision\footnote{Note that $\bx\bx^T$ is a symmetric matrix of quadratic monomials.}, i.e.
\begin{align}
\sum_{a=1}^n \sum_{b=1}^n \mu_{ab} &= 1, & \sum_{a=1}^n \sum_{b=1}^n\bv_a\mu_{ab} &= \bx & {\rm and} & & \sum_{a=1}^n \sum_{b=1}^n\bv_a\bv_b^T\mu_{ab} = \bx \bx^T.\label{eq:cplpqp}
\end{align}
\end{proposition}
\begin{proof}
The result is immediate from properties B\ref{b:partition} and B\ref{b:linprec} of the $\lambda_i$ functions.
\end{proof}

The product rule ensures that Theorem~\ref{th:aicmresult} extends immediately to the pairwise product functions.

\begin{corollary}\label{cr:aicmcor}
Let $\Omega$ be a convex polygon satisfying G1, G2, and G3, and let $\lambda_i$ denote a set of barycentric coordinates satisfying the result of Theorem~\ref{th:aicmresult} (e.g.\  $\ltri$, $\lopt$, $\lwach$, $\lsib$, or $\lmval$).
Then pairwise products of the $\lambda_i$ functions are all bounded in $H^1$, i.e. there exists a constant $C>0$ such that
\begin{equation}\label{eq:aicmcor}
\hpn{\mu_{ab}}{1}{\Omega} \leq C.
\end{equation}
\end{corollary}

While the $\{\mu_{ab}\}$ functions are commonly used on triangles to provide a quadratic Lagrange element, they have not been considered in the context of generalized barycentric coordinates on convex polygons as considered here.
Langer and Seidel have considered higher order barycentric interpolation in the computer graphics literature~\cite{LS08};
their approach, however, is for problems requiring $C^1$-continuous interpolation rather than the weaker $H^1$-continuity required for finite element theory.

In the remainder of this section, we describe notation that will be used to index functions throughout the rest of the paper.
Since $\mu_{ab}=\mu_{ba}$, the summations from (\ref{eq:cplpqp}) can be written in a symmetric expansion.
Define the paired index set
\[
I := \left\{ \{a,b\} \,|\, a,b\in \{1, \ldots, n\} \right\}.
\]
Note that sets with cardinality $1$ occur when $a=b$ and \textit{are} included in $I$.
We partition $I$ into three subsets corresponding to geometrical features of the polygon: vertices, edges of the boundary, and interior diagonals.
More precisely, $I= V \cup E \cup D$, a disjoint union, where
\begin{align*}
V & := \left\{ \{a,a\} \,|\, a\in \{1, \ldots, n\} \right\}; \\
E & := \left\{ \{a,a+1\} \,|\, a\in \{1, \ldots, n\} \right\};\\
D & := I \setminus \left(V\cup E \right).
\end{align*}
In the definition of $E$ above (and in general for indices throughout the paper), values are interpreted modulo $n$, i.e.\ $\{n,n+1\}$,  $\{n,1\}$, and $\{0,1\}$ all correspond to the edge between vertex $1$ and vertex $n$.
To simplify notation, we will omit the braces and commas when referring to elements of the index set $I$.
For instance, instead of $\mu_{\{a,b\}}$, we write just $\mu_{ab}$.
We emphasize that $ab \in I$ refers to an unordered and possibly non-distinct pair of vertices. 
Occasionally we will also use the abbreviated notation
\[\bv_{ab}:= \frac{\bv_a+\bv_b}{2},\]
so that $\bv_{aa}$ is just a different expression for $\bv_a$.
Under these conventions, the precision properties from (\ref{eq:cplpqp}) can be rewritten as follows.\newpage
\renewcommand{\labelenumi}{Q\arabic{enumi}.}
\begin{enumerate}
\item \textbf{Constant Precision}:   $\ds\sum_{aa\in V}\mu_{aa}+\sum_{ab\in E\cup D}2\mu_{ab}=1$ \label{q:pou}
\item \textbf{Linear Precision}:  $\ds\sum_{aa\in V}\bv_{aa}\mu_{aa}+\sum_{ab\in E\cup D}2\bv_{ab}\mu_{ab}=\bx$  \label{q:lp}
\item \textbf{Quadratic Precision}:  $\ds\sum_{aa\in V}\bv_a\bv_a^T\mu_{aa}+\sum_{ab\in E\cup D}(\bv_a\bv_b^T+\bv_b\bv_a^T)\mu_{ab}=\bx\bx^T$ \label{q:qp}
\end{enumerate}

\section{Reducing Quadratic Elements to Serendipity Elements}
\label{sec:construct}

We now seek to reduce the set of pairwise product functions $\{\mu_{ab}\}$ to a basis $\{\xi_{ij}\}$ for a serendipity quadratic finite element space.
Our desired basis must
\renewcommand{\labelenumi}{(\roman{enumi})}
\begin{enumerate}
\item span all quadratic polynomials of two variables on $\Omega$,
\item be exactly the space of quadratic polynomials (of one variable) when restricted to edges of $\Omega$, and
\item contain only $2n$ basis functions.
\end{enumerate}
The intuition for how to achieve this is seen from the number of distinct pairwise products:
\[|\{\mu_{ab}\}| = |I| = |V|+|E|+|D| = n + n + \frac {n(n-3)}2 = n +{n\choose 2}\]
On $\p\Omega$, functions with indices in $V$ vanish on all but two adjacent edges, functions with indices in $E$ vanish on all but one edge, and functions with indices in $D$ vanish on all edges.
Since Q1-Q3 hold on all of $\Omega$, including $\p\Omega$, the set $\{\mu_{ab} : ab\in V\cup E\}$ satisfies (ii) and (iii), but not necessarily (i).
Thus, our goal is to add linear combinations of functions with indices in $D$ to those with indices in $V$ or $E$ such that (i) is ensured.

We formalize this goal as a linear algebra problem:\ find a matrix $\A$ for the equation
\begin{equation}
\label{eq:Abeqg}
[\xi_{ij}] := \A [\mu_{ab}]
\end{equation}
such that $[\xi_{ij}]$ satisfies the following conditions analogous to Q1-Q3:
\renewcommand{\labelenumi}{$\Qx$\arabic{enumi}.}
\begin{enumerate}
\item \textbf{Constant Precision}:  $\ds\sum_{ii\in V} \xi_{ii}+\sum_{i(i+1)\in E}2\xi_{i(i+1)}=1$. \label{qg:pou}
\item \textbf{Linear Precision}: $\ds\sum_{ii\in V}\bv_{ii}\xi_{ii}+\sum_{i(i+1)\in E} 2\bv_{i(i+1)}\xi_{i(i+1)}=\bx$.  \label{qg:lp}
\item \textbf{Quadratic Precision}:\\  $\ds\sum_{ii\in V}\bv_i\bv_i^T\xi_{ii}+\sum_{i(i+1)\in E}(\bv_i\bv_{i+1}^T+\bv_{i+1}\bv_i^T)\xi_{i(i+1)}=\bx\bx^T$. \label{qg:qp}
\end{enumerate}
Since (\ref{eq:Abeqg}) is a linear relationship, we are still able to restrict our analysis to a reference set of unit diameter polygons (recall Figure~\ref{fg:affinetrans}). 
Specifically if matrix $\A$ yields a ``reference'' basis $T[\xi_{ij}]=\A T[\mu_{ab}]$ satisfying $\Qx$1-$\Qx$3, then the ``physical'' basis $[\xi_{ij}]=\A [\mu_{ab}]$ also satisfies $\Qx$1-$\Qx$3.

To specify $\A$ in (\ref{eq:Abeqg}), we will use the specific basis orderings
\begin{align}
[\xi_{ij}] &:= [~\underbrace{\xi_{11}, \xi_{22},\ldots,\xi_{nn}}_{\text{indices in $V$}},~~\underbrace{\xi_{12},\xi_{23},\ldots,\xi_{(n-1)n},\xi_{n(n+1)}}_{\text{indices in $E$}}~],\label{eq:order-xi}\\
[\mu_{ab}] & := [~\underbrace{\mu_{11}, \mu_{22},\ldots,\mu_{nn}}_{\text{indices in $V$}},~~\underbrace{\mu_{12},\mu_{23},\ldots,\mu_{(n-1)n},\mu_{n(n+1)}}_{\text{indices in $E$}}, \label{eq:order-mu}\\
& \qquad\qquad\qquad\qquad \underbrace{\mu_{13},\ldots, \text{(lexicographical)},\ldots,\mu_{(n-2)n}}_{\text{indices in $D$}}~]. \notag
\end{align}
The entries of $\A$ are denoted $c^{ij}_{ab}$ following the orderings given in (\ref{eq:order-xi})-(\ref{eq:order-mu}) so that
\begin{equation}
\label{eq:A-entries}
\A:= \left[ 
\begin{array}{ccccc}
c^{11}_{11} & \cdots & c^{11}_{ab} & \cdots &  c^{11}_{(n-2)n} \\ 
\vdots & \ddots & \vdots & \ddots  & \vdots \\
c^{ij}_{11} & \cdots & c^{ij}_{ab} & \cdots &  c^{ij}_{(n-2)n} \\ 
\vdots & \ddots & \vdots & \ddots  & \vdots \\
c^{n(n+1)}_{11} & \cdots & c^{n(n+1)}_{ab} & \cdots &  c^{n(n+1)}_{(n-2)n}
\end{array}
\right].
\end{equation}
A sufficient set of constraints on the coefficients of $\A$ to ensure $\Qx1$-$\Qx3$ is given by the following lemma.
\begin{lemma}
\label{lem:qc-implies-qx}
The constraints $\Qc1$-$\Qc3$ listed below imply $\Qx1$-$\Qx3$, respectively.  
That is, $\Qc 1 \Rightarrow \Qx 1$, $\Qc 2 \Rightarrow \Qx 2$, and $\Qc 3 \Rightarrow \Qx 3$. 
\renewcommand{\labelenumi}{$\Qc$\arabic{enumi}.}
\begin{enumerate}
\item $\ds \sum_{ii\in V}c^{ii}_{aa}+\sum_{i(i+1)\in E}2 c^{i(i+1)}_{aa}  = 1~\text{$\forall aa\in V$, and}\\
\sum_{ii\in V}c^{ii}_{ab}+\sum_{i(i+1)\in E}2 c^{i(i+1)}_{ab}  = 2,~\text{$\forall ab\in E\cup D$}.$ \label{qc:pou}
\item $\ds \sum_{ii\in V} c^{ii}_{aa}\bv_{ii}+\sum_{i(i+1)\in E}2 c^{i(i+1)}_{aa}\bv_{i(i+1)} = \bv_{aa}~\text{$\forall aa\in V$, and}\\ 
\sum_{ii\in V}c^{ii}_{ab}\bv_{ii}+\sum_{i(i+1)\in E}2 c^{i(i+1)}_{ab}\bv_{i(i+1)}  = 2\bv_{ab},~\text{$\forall ab\in E\cup D$}.$ \label{qc:lp}
\item $\ds \sum_{ii\in V}c^{ii}_{aa}\bv_i\bv_i^T+ \sum_{i(i+1)\in E}c^{i(i+1)}_{aa}(\bv_i\bv_{i+1}^T+\bv_{i+1}\bv_i^T)  = \bv_a\bv_a^T~\text{$\forall a\in V$, and}\\ 
\sum_{ii\in V}c^{ii}_{ab}\bv_i\bv_i^T +\sum_{i(i+1)\in E}c^{i(i+1)}_{ab}(\bv_i\bv_{i+1}^T+\bv_{i+1}\bv_i^T)  = \bv_a\bv_b^T+\bv_b\bv_a^T,~\text{$\forall ab\in E\cup D$}.$ \label{qg:qp}
\end{enumerate}
\end{lemma}
\begin{proof}
Suppose $\Qc1$ holds.
Substituting the expressions from $\Qc1$ into the coefficients of Q1 (from the end of Section~\ref{sec:bkgdnot}), we get
\begin{equation*}\begin{split}
\sum_{aa\in V}\left(\sum_{ii\in V} c^{ii}_{aa}+
 \sum_{i(i+1)\in E}2 c^{i(i+1)}_{aa} \right) \mu_{aa}+ & \\
\sum_{ab\in E\cup D}\left(\sum_{ii\in V} c^{ii}_{ab}+ \sum_{i(i+1)\in E}2 c^{i(i+1)}_{ab}\right) \mu_{ab} & =1.
\end{split}
\end{equation*}
Regrouping this summation over $ij$ indices instead of $ab$ indices, we have
\begin{equation}
\label{eq:pou-q1-exp}
\sum_{ii\in V}\left(\sum_{ab\in I}c^{ii}_{ab}\mu_{ab}\right) + \sum_{i(i+1)\in E} 2\left(\sum_{ab\in I}c^{i(i+1)}_{ab}\mu_{ab}\right)=1.
\end{equation}
Since (\ref{eq:Abeqg}) defines $\ds \xi_{ij} = \sum_{ab\in I}c^{ij}_{ab}\mu_{ab}$, (\ref{eq:pou-q1-exp}) is exactly the statement of $\Qx1$.
The other two cases follow by the same technique of regrouping summations.
\end{proof}
We now give some remarks about our approach to finding coefficients satisfying $\Qc$1-$\Qc$3.
Observe that the first equation in each of $\Qc$1-$\Qc$3 is satisfied by
\begin{equation}
\label{eq:A-struc-1}
c^{ii}_{aa}:=\delta_{ia}\quad\text{and}\quad c^{i(i+1)}_{aa}:= 0
\end{equation}
Further, if $ab=a(a+1)\in E$, the second equation in each of $\Qc$1-$\Qc$3 is satisfied by
\begin{equation}
\label{eq:A-struc-2}
c^{ii}_{a(a+1)}:=0\quad\text{and}\quad c^{i(i+1)}_{a(a+1)}:= \delta_{ia}
\end{equation}
The choices in (\ref{eq:A-struc-1}) and (\ref{eq:A-struc-2}) give $\A$ the simple structure
\begin{equation}
\label{eq:A-struc}
\A:= \left[ 
\begin{array}{c|c} \id & \A' \end{array}
\right],
\end{equation}
where $\id$ is the $2n\times 2n$ identity matrix.
Note that this corresponds exactly to our intuitive approach of setting each $\xi_{ij}$ function to be the corresponding $\mu_{ij}$ function plus a linear combination of $\mu_{ab}$ functions with $ab\in D$.
Also, with this selection, we can verify that many of the conditions which are part of $\Qc$1, $\Qc$2 and $\Qc$3 hold. 
Specifically, whenever $ab\in V\cup E$, the corresponding conditions hold as stated in the following lemma.

\begin{lemma}\label{lem:easycases}
The first $2n$ columns of the matrix $\A$ given by (\ref{eq:A-struc}), i.e., the identity portion, ensure $\Qc$1, $\Qc$2 and $\Qc$3 hold for $ab\in V\cup E$.
\end{lemma}
\begin{proof}
This lemma follows from direct substitution. In each case, there is only one nonzero element $c_{aa}^{aa}$ or $c_{a(a+1)}^{a(a+1)}$ on the hand side of the equation from $\Qc$1, $\Qc$2 or $\Qc$3 and substituting $1$ for that coefficient gives the desired equality.
\end{proof}

\begin{figure}
\centering
\includegraphics[width=.4\textwidth]{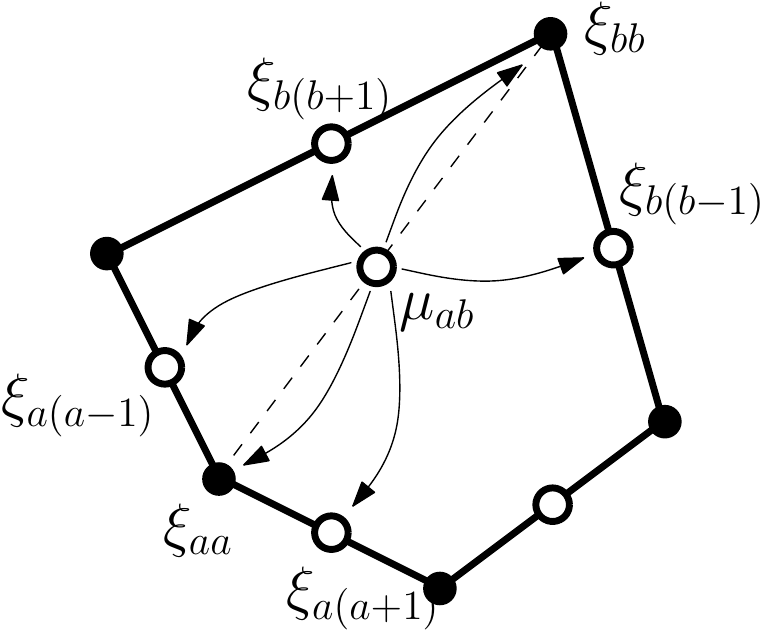}
\caption{When constructing the matrix $\A$, only six non-zero elements are used in each column corresponding to an interior diagonal of the pairwise product basis. In the serendipity basis, the interior diagonal function $\mu_{ab}$ only contributes to six basis functions as shown, corresponding to the vertices of the diagonal's endpoints and the midpoints of adjacent boundary edges. }\label{fg:selectingnonzeros}
\end{figure}

It remains to define $\A'$, i.e.\ those coefficients $c_{ab}^{ij}$ with $ab \in D$ and verify the corresponding equations in $\Qc$1, $\Qc$2, and $\Qc$3.
For each column of $\A'$, $\Qc$1, $\Qc$2, and $\Qc$3 yield a system of six scalar equations for the $2n$ variables $\{ c_{ab}^{ij} \}_{ij\in V\cup E}$.
Since we have many more variables than equations, there remains significant flexibility in the construction of a solution.
In the upcoming sections, we will present such a solution where all but six of the coefficients in each column of $\A'$ are set to zero.
The non-zero coefficients are chosen to be $c_{ab}^{a(a-1)}$, $c_{ab}^{aa}$, $c_{ab}^{a(a+1)}$, $c_{ab}^{b(b-1)}$, $c_{ab}^{bb}$, and $c_{ab}^{b(b+1)}$ as these have a natural correspondence to the geometry of the polygon and the edge $ab$; see Figure~\ref{fg:selectingnonzeros}.

We will show that the system of equations $\Qc$1-$\Qc$3 with this selection of non-zero coefficients for $\A'$ has an explicitly constructible solution.
The solution is presented for special classes of polygons in Section~\ref{sec:specialcase} and for generic convex polygons in Section~\ref{sec:genpoly}.
In each case, we prove a uniform bound on the size of the coefficients of $\A$, a sufficient result to control $\hpn{\xi_{ij}}{1}{\Omega}$, as the following lemma shows.

\begin{lemma}
\label{lem:cbd-implies-Abd}
Let $\Omega$ be a convex polygon satisfying G1, G2, and G3, and let $\lambda_i$ denote a set of barycentric coordinates satisfying the result of Theorem~\ref{th:aicmresult} (e.g.\  $\ltri$, $\lopt$, $\lwach$, $\lsib$, or $\lmval$).
Suppose there exists $M>1$ such that for all entries of $\A'$, $|c^{ij}_{ab}|<M$.  
Then the functions $\xi_{ij}$ are all bounded in $H^1$, i.e. there exists a constant $B>0$ such that
\begin{equation}
\label{eq:xi-bd}
\hpn{\xi_{ij}}{1}{\Omega} \leq B.
\end{equation}
\end{lemma}
\begin{proof}
Since $\xi_{ij}$ is defined by (\ref{eq:Abeqg}), Corollary~\ref{cr:aicmcor} implies that there exists $C>0$ such that
\[\hpn{\xi_{ij}}{1}{\Omega} \leq ||\A|| \max_{ab}\hpn{\mu_{ab}}{1}{\Omega}< C||\A||.\]
Since the space of linear transformations from $\R^{n(n+1)/2}$ to $\R^{2n}$ is finite-dimensional, all norms on $\A$ are equivalent.
Thus, without loss of generality, we interpret $||\A||$ as the maximum absolute row sum norm, i.e.
\begin{equation}
\label{eq:Anorm}
||\A||:= \max_{ij}\sum_{ab}|c_{ab}^{ij}|.
\end{equation}
By the structure of $\A$ from (\ref{eq:A-struc}) and the hypothesis, we have 
\[||\A||\leq \frac{n(n+1)}{2}M\]
\end{proof}

\section{Special Cases of the Serendipity Reduction}
\label{sec:specialcase}

Before showing that $\Qc$1-$\Qc$3 can be satisfied in a general setting, we study some simpler special cases in which symmetry reduces the number of equations that must be satisfied simultaneously.

\subsection{Unit Square}
We begin with the case where serendipity elements were first examined, namely over meshes of squares.
In recent work by Arnold and Awanou~\cite{AA10}, the quadratic serendipity space on the unit square, denoted $\cS_2(I^2)$, is defined as the span of eight monomials:
\begin{equation}
\label{eq:def-s2}
\cS_2(I^2) := \text{span}\left\{1, x, y, x^2, xy, y^2, x^2y, xy^2\right\}
\end{equation}
We will now show how our construction process recovers this same space of monomials.
Denote vertices on $[0,1]^2$ by
\[ \bv_1 = (0,0)\quad \bv_2 =(1,0) \quad  \bv_3= (1,1) \quad \bv_4=(0,1)\]
\begin{equation}
\label{eq:vtx-def-sqcase}
 \bv_{12} = (1/2,0) \quad \bv_{23} = (1,1/2) \quad   \bv_{34} = (1/2,1) \quad \bv_{14} = (0,1/2)
\end{equation}
\[ \bv_{13} = \bv_{24} = (1/2,1/2)\]
The standard bilinear basis for the square is 
\begin{align*}
\lambda_1 &= (1-x)(1-y)  & \lambda_2 &= x(1-y)\\
\lambda_4 &= (1-x)y  & \lambda_3 &= xy
\end{align*}
Since the $\lambda_i$ have vanishing second derivatives and satisfy the definition of barycentric coordinates, they are in fact the harmonic coordinates $\lopt$ in this special case.
Pairwise products give us the following $10$ (not linearly independent) functions
\begin{align*}
\mu_{11} & = (1-x)^2(1-y)^2 & \mu_{12} & = (1-x)x(1-y)^2\\ 
\mu_{22} & = x^2(1-y)^2 & \mu_{23} & = x^2(1-y)y\\ 
\mu_{33} & = x^2y^2 & \mu_{34} & = (1-x)x y^2\\ 
\mu_{44} & = (1-x)^2y^2 & \mu_{14} & = (1-x)^2(1-y)y\\ 
\mu_{13} & = (1-x)x(1-y)y & \mu_{24} & = (1-x)x(1-y)y
\end{align*}
For the \emph{special} geometry of the square, $\mu_{13} = \mu_{24}$, but this is not true for general quadrilaterals as we see in Section~\ref{subsec:quads}.
The serendipity construction eliminates the functions $\mu_{13}$ and $\mu_{24}$ to give an 8-dimensional space.
The basis reduction via the $\A$ matrix is given by
\begin{equation}
\label{eq:A-def-sqcase}
\left[\begin{array}{c}\xi_{11}\\ \xi_{22}\\\xi_{33}\\\xi_{44}\\\xi_{12}\\\xi_{23}\\\xi_{34}\\\xi_{14} \end{array} \right]
=
\left[
\begin{array}{cccccccccc}
1 & 0 & 0 & 0 & 0 & 0 & 0 & 0 & -1 & 0 \\ 
0 & 1 & 0 & 0 & 0 & 0 & 0 & 0 & 0 & -1 \\ 
0 & 0 & 1 & 0 & 0 & 0 & 0 & 0 & -1 & 0 \\ 
0 & 0 & 0 & 1 & 0 & 0 & 0 & 0 & 0 & -1 \\ 
0 & 0 & 0 & 0 & 1 & 0 & 0 & 0 & 1/2 & 1/2 \\ 
0 & 0 & 0 & 0 & 0 & 1 & 0 & 0 & 1/2 & 1/2 \\ 
0 & 0 & 0 & 0 & 0 & 0 & 1 & 0 & 1/2 & 1/2 \\ 
0 & 0 & 0 & 0 & 0 & 0 & 0 & 1 & 1/2 & 1/2 
\end{array}
\right]
\left[\begin{array}{c}\mu_{11}\\ \mu_{22}\\\mu_{33}\\\mu_{44}\\\mu_{12}\\\mu_{23}\\\mu_{34}\\\mu_{14}\\ \mu_{13}\\ \mu_{24} \end{array} \right]
\end{equation}
It can be confirmed directly that (\ref{eq:A-def-sqcase}) follows from the definitions of $\A$ given in the increasingly generic settings examined in Section~\ref{subsec:regpoly}, Section~\ref{subsec:quads} and Section~\ref{sec:genpoly}.
The resulting functions are
\begin{align}
\xi_{11} & = (1-x)(1-y)(1-x-y) & \xi_{12} & = (1-x)x(1-y) \label{eq:xi-def-sqcase} \\ 
\xi_{22} & = x(1-y)(x-y) & \xi_{23} & = x(1-y)y \notag \\ 
\xi_{33} & = xy(-1+x+y) & \xi_{34} & = (1-x)x y \notag \\ 
\xi_{44} & = (1-x)y(y-x) & \xi_{14} & = (1-x)(1-y)y \notag
\end{align}

\begin{theorem}
\label{thm:ser-unitsq}
For the unit square, the basis functions $\{\xi_{ij}\}$ defined in (\ref{eq:xi-def-sqcase}) satisfy $\Qx1$-$\Qx3$.
\end{theorem}
\begin{proof}
A simple proof is to observe that the coefficients of the matrix in (\ref{eq:A-def-sqcase}) satisfy $\Qc$1-$\Qc$3 and then apply Lemma~\ref{lem:qc-implies-qx}.
To illuminate the construction in this special case of common interest, we state some explicit calculations.
The constant precision condition $\Qx1$ is verified by the calculation
\begin{align*}
\xi_{11}+\xi_{22}+\xi_{33}+\xi_{44}+2\xi_{12}+2\xi_{23}+2\xi_{34}+2\xi_{14} = 1.
\end{align*}
The $x$ component of the linear precision condition $\Qx2$ is verified by the calculation
\begin{align*}
(\bv_{1})_x \xi_{11} +  (\bv_{2})_x \xi_{22} & + (\bv_{3})_x \xi_{33} + (\bv_{4})_x \xi_{44} +\\
2(\bv_{12})_x \xi_{12} &+ 2(\bv_{23})_x \xi_{23} + 2(\bv_{34})_x \xi_{34} + 2(\bv_{14})_x \xi_{14}\\
& = \xi_{22} + \xi_{33} + 2\cdot \frac{1}{2} \xi_{12} + 2\cdot 1\xi_{23} + 2\cdot \frac{1}{2} \xi_{34} \\
& = x.
\end{align*}
The verification for the $y$ component is similar. 
The $xy$ component of the quadratic precision condition $\Qx3$ is verified by
\begin{align*}
& (\bv_{1})_x(\bv_{1})_y\xi_{11} + 
(\bv_{2})_x(\bv_{2})_y\xi_{22} +
(\bv_{3})_x(\bv_{3})_y\xi_{33} +
(\bv_{4})_x(\bv_{4})_y\xi_{44}\\
& + \left[(\bv_{1})_x(\bv_{2})_y+(\bv_{2})_x(\bv_{1})_y\right]\xi_{12} +
\left[(\bv_{2})_x(\bv_{3})_y+(\bv_{3})_x(\bv_{2})_y\right]\xi_{23}\\
& + \left[(\bv_{3})_x(\bv_{4})_y+(\bv_{4})_x(\bv_{3})_y\right]\xi_{34} +
\left[(\bv_{4})_x(\bv_{1})_y+(\bv_{1})_x(\bv_{4})_y\right]\xi_{14}\\
& = \xi_{33} + \xi_{23} + \xi_{34} = xy.
\end{align*}
The monomials $x^2$ and $y^2$ can be expressed as a linear combination of the $\xi_{ij}$ similarly, via the formula given in $\Qx3$.
\end{proof}
\begin{corollary}
The span of the $\xi_{ij}$ functions defined by (\ref{eq:xi-def-sqcase}) is the standard serendipity space, i.e.\
\[{\rm span}\left\{\xi_{ii},\xi_{i(i+1)}\right\} =  \cS_2(I^2)\]
\end{corollary}
\begin{proof}
Observe that $x^2y = \xi_{23} + \xi_{33}$ and $xy^2 = \xi_{33} + \xi_{34}$.
By the definition of $\cS_2(I^2)$ in (\ref{eq:def-s2}) and the theorem, ${\rm span}\left\{\xi_{ii},\xi_{i(i+1)}\right\} \supset \cS_2(I^2)$.
Since both spaces are dimension eight, they are identical.
\end{proof}

\subsection{Regular Polygons}
\label{subsec:regpoly}

\begin{figure}
\centering
\includegraphics[width=.3\textwidth]{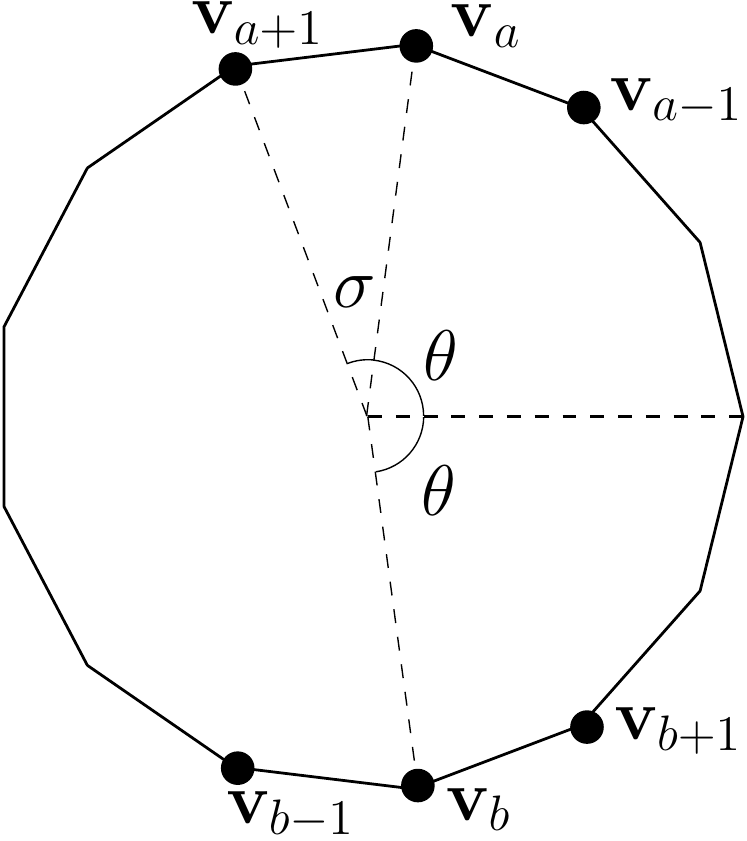}
\caption{Notation for the construction for a regular polygon.}\label{fg:regpoly}
\end{figure}

We now generalize our construction to any regular polygon with $n$ vertices.  
Without loss of generality, this configuration can be described by two parameters $0 < \sigma \leq \theta \leq \pi/2$ as shown in Figure~\ref{fg:regpoly}.  
Note that the $n$ vertices of the polygon are located at angles of the form $k\sigma$ where $k=0,1,\ldots,n-1$.  

For two generic non-adjacent vertices $\bv_a$ and $\bv_b$, the coordinates of the six relevant vertices (recalling Figure~\ref{fg:selectingnonzeros}) are:
\begin{align*}
\bv_a &=
\begin{bmatrix}
\cos\theta \\
\sin\theta
\end{bmatrix};
&
\bv_{a-1} &=
\begin{bmatrix}
\cos\tms \\
\sin\tms
\end{bmatrix};
&
\bv_{a+1} &=
\begin{bmatrix}
\cos\tps \\
\sin\tps
\end{bmatrix};\\
\bv_b &=
\begin{bmatrix}
\cos\theta \\
-\sin\theta
\end{bmatrix};
&
\bv_{b-1} &=
\begin{bmatrix}
\cos\tps \\
-\sin\tps
\end{bmatrix};
&
\bv_{b+1}&=
\begin{bmatrix}
\cos\tms \\
-\sin\tms
\end{bmatrix}.
\end{align*}

We seek to establish the existence of suitable constants $c^{aa}_{ab}$, $c^{a,a+1}_{ab}$, $c^{a-1,a}_{ab}$, $c^{bb}_{ab}$, $c^{b-1,b}_{ab}$, $c^{b,b+1}_{ab}$ which preserve quadratic precision and to investigate the geometric conditions under which these constants become large.
The symmetry of this configuration suggests that $c^{aa}_{ab}=c^{bb}_{ab}$, $c^{a-1,a}_{ab}=c^{b,b+1}_{ab}$, and $c^{a,a+1}_{ab}=c^{b-1,b}_{ab}$ are reasonable requirements.
For simplicity we will denote these constants by $c_0 := c^{aa}_{ab}$, $c_- := c^{a-1,a}_{ab}$, and $c_+ := c^{a,a+1}_{ab}$.

Thus equation $\Qc$1 (which contains only six non-zero elements) reduces to:
\begin{align}\label{eq:regpolyunity}
2c_0 + 4c_- + 4c_+ = 2.
\end{align}

$\Qc$2 involves two equations, one of which is trivially satisfied in our symmetric configuration.
Thus, the only restriction to maintain is
\begin{align}\label{eq:regpolylin}
2\cos\theta c_0 + 2\left[\cos\theta + \cos(\theta-\sigma)\right]c_- + 2\left[\cos\theta + \cos(\theta+\sigma)\right]c_+ = 2\cos\theta.
\end{align}

$\Qc$3 gives three more requirements, one of which is again trivially satisfied.  This gives two remaining restrictions:
\begin{align}
2\cos^2\theta c_0 + 4\cos\theta\cos\tms c_- + 4\cos\theta \cos\tps c_+ &= 2\cos^2\theta ; \label{eq:regpolyquad1}\\
2\sin^2\theta c_0 + 4\sin\theta\sin\tms c_- + 4\sin\theta \sin\tps c_+ &= -2\sin^2\theta . \label{eq:regpolyquad2}
\end{align}

Now we have four equations (\ref{eq:regpolyunity})-(\ref{eq:regpolyquad2}) and three unknowns $c_0$, $c_-$ and $c_+$.  Fortunately, equation (\ref{eq:regpolylin}) is a simple linear combination of (\ref{eq:regpolyunity}) and(\ref{eq:regpolyquad1}); specifically (\ref{eq:regpolylin}) is $\frac{\cos\theta}{2}$ times (\ref{eq:regpolyunity}) plus $\frac{1}{2\cos\theta}$ times (\ref{eq:regpolyquad1}).  With a little algebra, we can produce the system:
\begin{equation}
\left[
\begin{array}{ccc}
1 & 2 & 2\\
1 & 2(\cos\sigma+\sin\sigma\tan\theta) & 2(\cos\sigma-\sin\sigma\tan\theta)\\
1 & 2(\cos\sigma-\sin\sigma\cot\theta) & 2(\cos\sigma+\sin\sigma\cot\theta)
\end{array}
\right]
\left[
\begin{array}{c}
c_0\\
c_-\\
c_+
\end{array}
\right]
=
\left[
\begin{array}{c}
1\\
1\\
-1
\end{array}
\right].
\end{equation}
The solution of this system can be computed:
\begin{align*}
c_0 & = \frac{(-1+\cos\sigma) \cot\theta+(1+\cos\sigma) \tan\theta}{(-1+\cos\sigma) (\cot\theta+\tan\theta)};
\end{align*}
\vspace{-.18in}
\begin{align*}
c_- &= \frac{\cos\sigma - \sin\sigma \tan\theta - 1}{2\left(\tan\theta + \cot\theta\right)\sin\sigma\left(\cos\sigma-1\right)}; &
c_+ &= \frac{1-\cos\sigma - \sin\sigma \tan\theta}{2\left(\tan\theta + \cot\theta\right)\sin\sigma\left(\cos\sigma-1\right)}.
\end{align*}

Although $\tan\theta$ (and thus the solution above) is not defined for $\theta=\pi/2$, the solution in this boundary case can be defined by the limiting value which always exists.  We can now prove the following.

\begin{theorem}
For any regular polygon, the basis functions $\{\xi_{ij}\}$ constructed using the coefficients $c_{ab}^{aa} = c_{ab}^{bb} = c_0$, $c_{ab}^{a-1,a}=c_{ab}^{b,b+1}=c_-$, $c_{ab}^{a,a+1}=c_{ab}^{b-1,b}=c_+$ satisfy $\Qx$1-$\Qx$3.
\end{theorem}

\begin{proof}
The construction above ensures that the solution satisfies $\Qc$1, $\Qc$2, and $\Qc$3. 
\end{proof}

The serendipity element for regular polygons can be used for meshes consisting of only one regular polygon or a finite number of regular polygons. 
The former occurs only in meshes of triangles, squares and hexagons as these are the only regular polygons that can tile the plane. 
On the other hand, many tilings consisting of several regular polygons can be constructed using multiple regular polygons.
Examples include the snub square tiling (octagons and squares), the truncated hexagonal tiling (dodecahedra and triangles), the rhombitrihexagonal tiling (hexagons, squares, and triangles), and the truncated trihexagonal tiling (dodecagons, hexagons, and squares); see e.g.\ \cite{C1989}.
The construction process outlined above opens up the possibility of finite element methods applied over these types of mixed-geometry meshes, a mostly unexplored field.

\subsection{Generic Quadrilaterals}
\label{subsec:quads}

\begin{figure}
\centering
\includegraphics[width=.4\textwidth]{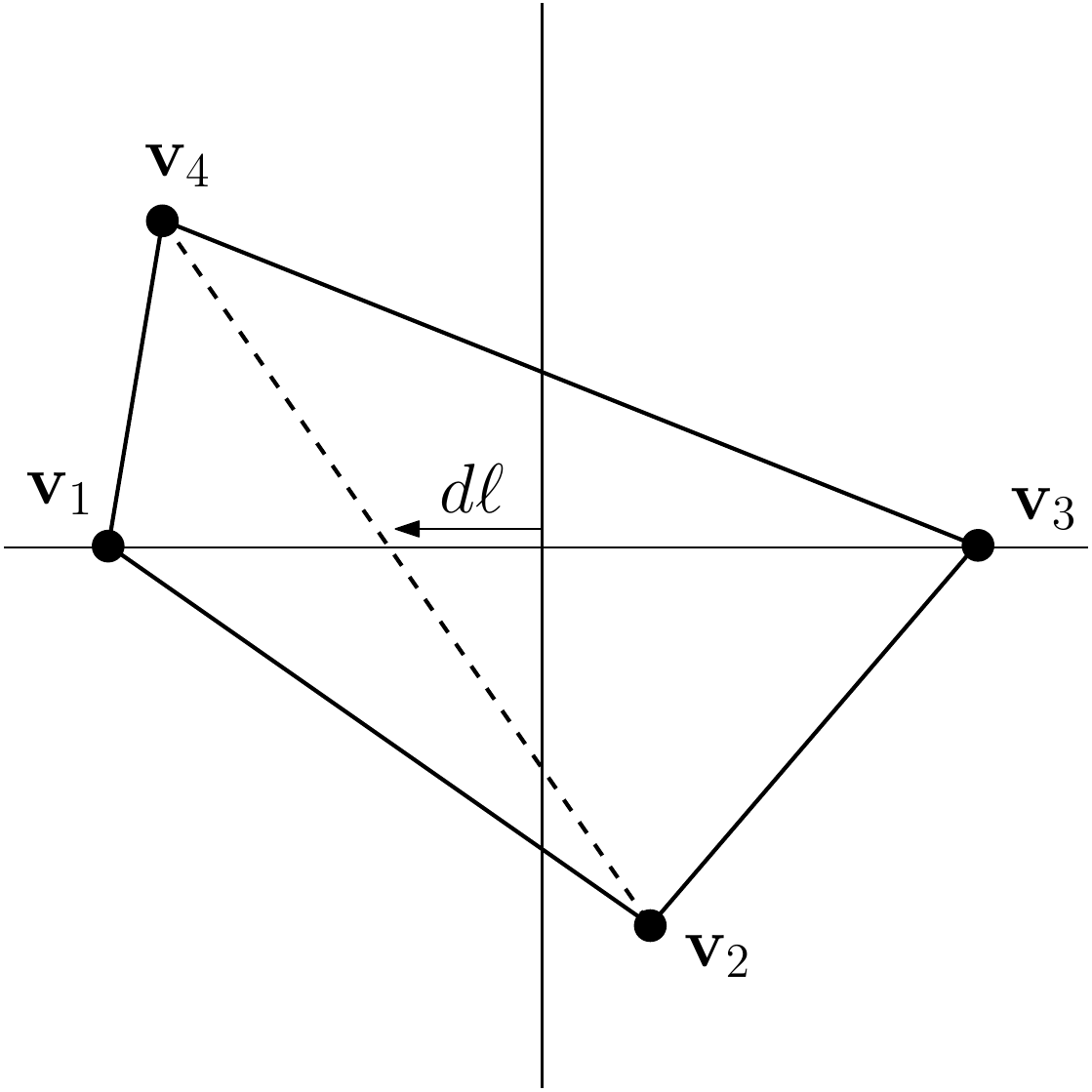}
\caption{A generic convex quadrilateral, rotated so that one of its diagonals lies on the $x$-axis. Geometrically, $c_{13}^{12}$ and $c_{13}^{34}$ are selected to be coefficients of the convex combination of $\bv_2$ and $\bv_4$ that lies on the $x$-axis.}
\label{fg:genquad}
\end{figure}

Fix a convex quadrilateral $\Omega$ with vertices $\bv_1$, $\bv_2$, $\bv_3$, and $\bv_4$, ordered counterclockwise. 
We will describe how to set the coefficients of the submatrix $\A'$ in (\ref{eq:A-struc}).
It suffices to describe how to set the coefficients in the `13'-column of the matrix, i.e., those of the form $c^{ij}_{13}$.  
The $`24'$-column can be filled using the same construction after permuting the indices.
Thus, without loss of generality, suppose that $\bv_1:=(-\ell,0)$ and $\bv_3:=(\ell,0)$ so that $\bv_2$ is below the $x$-axis and $\bv_4$ is above the $x$-axis, as shown in Figure~\ref{fg:genquad}.  
We have eight coefficients to set:
\[c^{11}_{13},\; c^{22}_{13},\; c^{33}_{13},\; c^{44}_{13},\;c^{12}_{13},\; c^{23}_{13},\; c^{34}_{13},\;{\rm and} \; c^{14}_{13}.\]
Using a subscript $x$ or $y$ to denote the corresponding component of a vertex, define the coefficients as follows.
\begin{alignat}{2}
\ds c^{22}_{13} & :=0 & \quad c^{44}_{13} & :=0 \label{eq:qcoef1}\\
\ds c^{12}_{13} & := \frac{(\bv_4)_y}{(\bv_4)_y-(\bv_2)_y} &  c^{34}_{13} & := \frac{(\bv_2)_y}{(\bv_2)_y-(\bv_4)_y} \label{eq:qcoef2}\\
\ds c^{23}_{13} & := c^{12}_{13} & c^{14}_{13} &:= c^{34}_{13} \label{eq:qcoef3}\\
\ds c^{11}_{13} & := \frac{c^{12}_{13}(\bv_2)_x + c^{34}_{13}(\bv_4)_x}{\ell}-1 &  c^{33}_{13} & := -\frac{c^{12}_{13}(\bv_2)_x + c^{34}_{13}(\bv_4)_x}{\ell}-1 \label{eq:qcoef4}
\end{alignat}
Note that there following the strategy shown in Figure~\ref{fg:selectingnonzeros}, there are only six non-zero entries.
For ease of notation in the rest of this section, we define the quantity
\[d:=\frac{c^{12}_{13}(\bv_2)_x+c^{34}_{13}(\bv_4)_x}{\ell}.  \]
First we assert that the resulting basis does span all quadratic polynomials.

\begin{theorem}
For any quadrilateral, the basis functions $\{\xi_{ij}\}$ constructed using the coefficients given in (\ref{eq:qcoef1})-(\ref{eq:qcoef4}) satisfy $\Qx$1-$\Qx$3.
\end{theorem}
\begin{proof}
Considering Lemmas~\ref{lem:qc-implies-qx} and \ref{lem:easycases}, we only must verify $\Qc$1-$\Qc$3 in the cases when $ab\in D = \{13,24\}$. 
This will be verified directly by substituting (\ref{eq:qcoef1})-(\ref{eq:qcoef4}) into the constraints $\Qc$1-$\Qc$3 in the case $ab=13$.  
As noted before, the $ab=24$ case is identical, requiring only a permutation of indices.
First note that
\begin{equation}
\label{eq:c11c33con}
c^{11}_{13}+c^{33}_{13}=-2\quad\text{and}\quad  c^{11}_{13}-c^{33}_{13}= 2d.
\end{equation}
For $\Qc$1, the sum reduces to
\[c^{11}_{13}+c^{33}_{13}+4(c^{12}_{13}+c^{34}_{13})=-2+4(1)=2,\]
as required.  For $\Qc$2, the $x$-coordinate equation reduces to
\[\ell(c^{33}_{13}-c^{11}_{13})+2d\ell = 0\]
by (\ref{eq:c11c33con}) which is the desired inequality since we fixed (without loss of generality) $\bv_ab = (0,0)$.   The $y$-coordinate equation reduces to $2(c^{12}_{13}(\bv_2)_y+c^{34}_{13}(\bv_4)_y)=0$ which holds by (\ref{eq:qcoef2}).  Finally, a bit of algebra reduces the matrix equality of $\Qc$3 to only the equality $\ell^2(c^{11}_{13}+c^{33}_{13})=-2\ell^2$ of its first entry (all other entries are zero), which holds by (\ref{eq:c11c33con}).
\end{proof}

\begin{theorem}
Over all convex quadrilaterals, $||\A||$ is uniformly bounded.
\end{theorem}
\begin{proof}
By Lemma~\ref{lem:cbd-implies-Abd}, it suffices to bound $|c^{ij}_{13}|$ uniformly.
First observe that the convex combination of the vertices $\bv_2$ and $\bv_4$ using coefficients $c^{12}_{13}$ and $c^{34}_{13}$ produces a point lying on the $x$-axis, i.e.,
\begin{align}
1 = & c^{12}_{13}+c^{34}_{13}, \,\, {\rm and} \label{eq:qcnvx1}  \\
0 = & c^{12}_{13}(\bv_2)_y+c^{34}_{13}(\bv_4)_y. \label{eq:qcnvx2}
\end{align}
Since $(\bv_2)_y>0$ and $(\bv_4)_y<0$, (\ref{eq:qcoef3}) implies that $c^{12}_{13},c^{34}_{13}\in (0,1)$.  By (\ref{eq:qcoef3}), it also follows that $ c^{23}_{13}, c^{14}_{13}\in (0,1)$.

For $c^{11}_{13}$ and $c^{33}_{13}$, note that the quantity $d\ell$ is the $x$-intercept of the line segment connecting $\bv_2$ and $\bv_4$.  Thus $d\ell\in[-\ell,\ell]$ by convexity.  So $d\in[-1,1]$ and thus (\ref{eq:qcoef4}) implies $|c^{11}_{13}|=|d-1|\leq 2$ and $|c^{33}_{13}|=|-d-1|\leq 2$.
\end{proof}

\section{Proof of the Serendipity Reduction on Generic Convex Polygons}
\label{sec:genpoly}

\begin{figure}
\centering
\includegraphics[width=.5\textwidth]{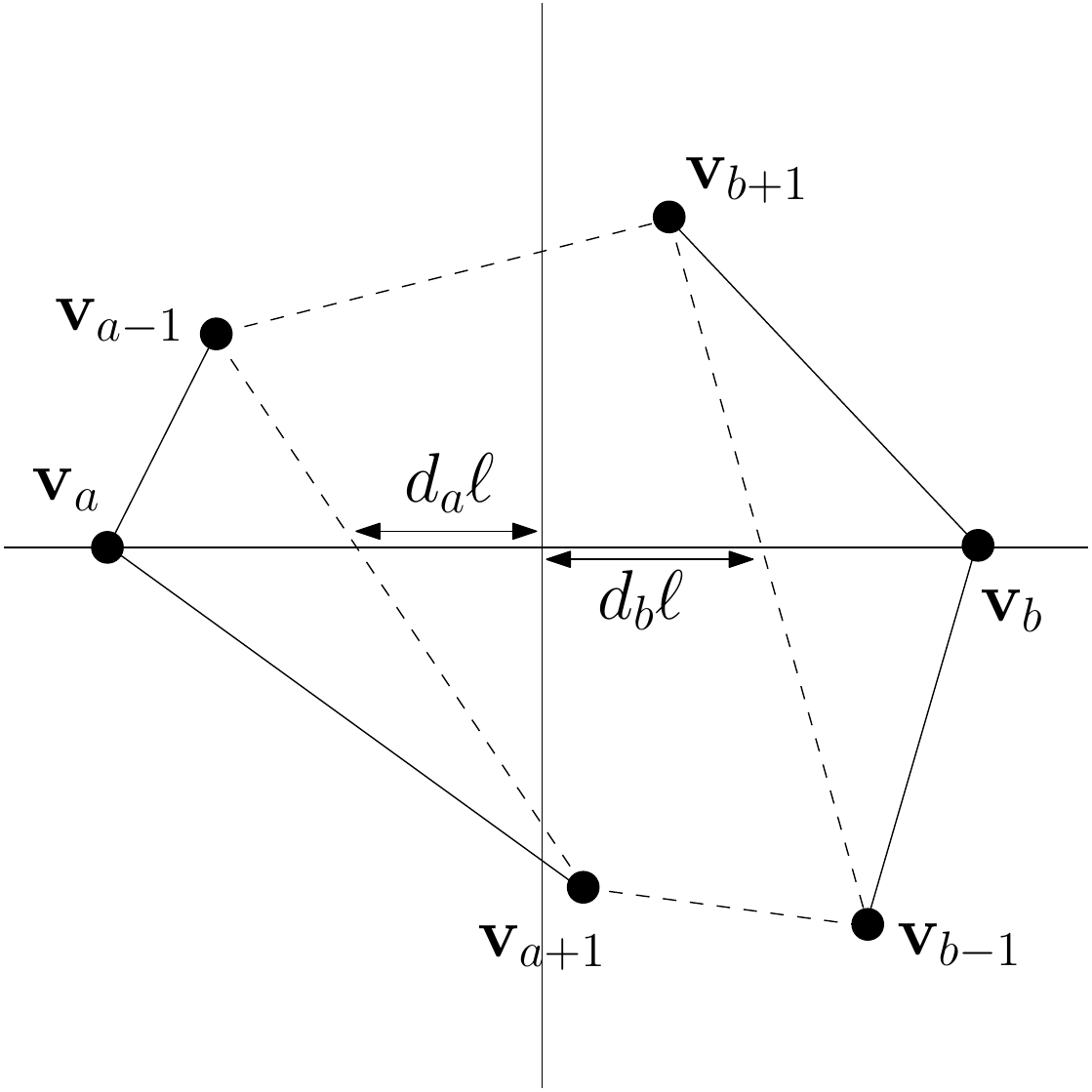}
\caption{Generic convex polygon, rotated so that $\bv_a=(-\ell,0)$ and $\bv_b=(\ell,0)$.  The $x$-intercept of the line between $\bv_{a-1}$ and $\bv_{a+1}$ is defined to be $-d_a\ell$ and the $x$-intercept of the line between $\bv_{b-1}$ and $\bv_{b+1}$ is defined to be $d_b\ell$.}
\label{fg:genpoly}
\end{figure}

We now define the sub-matrix $\A'$ from (\ref{eq:A-struc}) in the case of a generic polygon.  Pick a column of $\A'$, i.e., fix $ab \in D$.  The coefficients $c^{ij}_{ab}$ are constrained by a total of six equations $\Qc$1, $\Qc$2, and $\Qc$3.  
As before (recall Figure~\ref{fg:selectingnonzeros}), six non-zero coefficients will be selected in each column to satisfy these constraints. 
Specifically,
\begin{equation}
\label{eq:pre-set-coeffs}
c^{ii}_{ab}:=0,\;\text{for $i\not\in\{a,b\}$}\quad\text{and}\quad c^{i(i+1)}_{ab}=0,\;\text{for $i\not\in\{a-1,a,b-1,b\}$},
\end{equation}
leaving only the following six coefficients to be determined:
\[c^{aa}_{ab},\; c^{bb}_{ab},\; c^{(a-1)a}_{ab},\; c^{a(a+1)}_{ab},\; c^{(b-1)b}_{ab},\; {\rm and} \; c^{b(b+1)}_{ab}.\]
For the remainder of this section, we will omit the subscript $ab$ to ease the notation.  Writing out $\Qc$1-$\Qc$3 for this fixed $ab$ pair, we have six equations with six unknowns:
\begin{align*}
\ds c^{aa}+c^{bb}+2c^{(a-1)a}+2c^{a(a+1)}+2c^{(b-1)b}+2c^{b(b+1)} & = 2; \\
\ds c^{aa}\bv_{aa} +2c^{(a-1)a}\bv_{(a-1)a} +2c^{a(a+1)}\bv_{a(a+1)} + & \\
\ds c^{bb}\bv_{bb} + 2 c^{(b-1)b}\bv_{(b-1)b}  +2c^{b(b+1)}\bv_{b(b+1)} & = 2\bv_{ab}; \\
\ds c^{aa}\bv_a\bv_a^T + c^{(a-1)a}(\bv_{a-1}\bv_a^T+\bv_a\bv_{a-1}^T) + c^{a(a+1)}(\bv_a\bv_{a+1}^T+\bv_{a+1}\bv_a^T)+  & \\
\ds c^{bb}\bv_b\bv_b^T + c^{(b-1)b}(\bv_{b-1}\bv_b^T+\bv_b\bv_{b-1}^T)  + c^{b(b+1)}(\bv_b\bv_{b+1}^T+\bv_{b+1}\bv_b^T) & = \bv_a\bv_b^T+\bv_b\bv_a^T.
\end{align*}

Assume without loss of generality that $\bv_a = (-\ell,0)$ and $\bv_b = (\ell,0)$ with $\ell < 1/2$ (since $\Omega$ has diameter $1$).  We introduce the terms $d_a$ and $d_b$ defined by
\begin{align}
d_a &:= \frac{(\bv_{a-1})_x(\bv_{a+1})_y-(\bv_{a+1})_x(\bv_{a-1})_y}{(\bv_{a-1})_y-(\bv_{a+1})_y}\cdot\frac 1{\ell}, \,\, {\rm and} \label{eq:dadef} \\
d_b &:= \frac{(\bv_{b+1})_x(\bv_{b-1})_y-(\bv_{b-1})_x(\bv_{b+1})_y}{(\bv_{b-1})_y-(\bv_{b+1})_y}\cdot\frac 1{\ell}. \label{eq:dbdef}
\end{align}
These terms have a concrete geometrical interpretation as shown in Figure~\ref{fg:genpoly}: $-d_a\ell$ is the $x$-intercept of the line between $\bv_{a-1}$ and $\bv_{a+1}$, while $d_b\ell$ is the $x$-intercept of the line between $\bv_{b-1}$ and $\bv_{b+1}$.  Thus, by the convexity assumption, $d_a,d_b\in[-1,1]$.  Additionally, $-d_a\leq d_b$ with equality only in the case of a quadrilateral which was dealt with previously.  For ease of notation and subsequent explanation, we also define
\begin{equation}
\label{eq:defofs}
s := \ds\frac {2}{2-(d_a+d_b)}.
\end{equation}
First we choose $c^{(a-1)a}$ and $c^{a(a+1)}$ as the solution to the following system of equations:
\begin{align}
c^{(a-1)a} + c^{a(a+1)} &= s; \label{eq:gen1}\\
c^{(a-1)a}\bv_{a-1} + c^{a(a+1)}\bv_{a+1} &= s d_a \bv_a. \label{eq:gen2}
\end{align}
There are a total of three equations since (\ref{eq:gen2}) equates vectors, but it can be verified directly that this system of equations is only rank two.  Moreover, any two of the equations from (\ref{eq:gen1}) and (\ref{eq:gen2}) suffice to give the same unique solution for $c^{(a-1)a}$ and $c^{a(a+1)}$.

Similarly, we select $c^{(b-1)b}$ and $c^{b(b-1)}$ as the solution to the system:
\begin{align}
c^{(b-1)b} + c^{b(b+1)} &= s; \label{eq:gen3}\\
c^{(b-1)b}\bv_{b-1} + c^{b(b+1)}\bv_{b+1} &= s d_b \bv_b. \label{eq:gen4}
\end{align}
Finally, we assign $c^{aa}$ and $c^{bb}$ by
\begin{align}
c^{aa} & = \ds\frac{-2-2d_a}{2-(d_a+d_b)} \,\,\, {\rm and} \label{eq:gen5} \\
c^{bb} & = \ds\frac{-2-2d_b}{2-(d_a+d_b)}, \label{eq:gen6}
\end{align}
and claim that this set of coefficients leads to a basis with quadratic precision.
\begin{theorem}
For any convex polygon, the basis functions $\{\xi_{ij}\}$ constructed using the coefficients defined by (\ref{eq:gen1})-(\ref{eq:gen6}) satisfy $\Qx$1-$\Qx$3.
\end{theorem}
\begin{proof}
Based on Lemmas~\ref{lem:qc-implies-qx} and \ref{lem:easycases}, it only remains to verify that $\Qc$1, $\Qc$2, and $\Qc$3 hold when $ab\in D$. 
Observe that $c^{aa}$ and $c^{bb}$ satisfy the following equations:
\begin{align}
c^{aa} + c^{bb} + 4s & = 2; \label{eq:gen7} \\
c^{aa} - c^{bb} + s(d_a-d_b) &= 0;  \label{eq:gen8} \\
c^{aa} + c^{bb} +2s(d_a+d_b) & = -2. \label{eq:gen9}
\end{align}
First, note that $\Qc$1 follows immediately from (\ref{eq:gen1}), (\ref{eq:gen3}) and (\ref{eq:gen7}).

The linear precision conditions ($\Qc2$) are just a matter of algebra.  Equations (\ref{eq:gen1})-(\ref{eq:gen4}) yield
\flushleft{
$c^{aa}\bv_{aa} +c^{bb}\bv_{bb}  +2c^{(a-1)a}\bv_{(a-1)a} +2c^{a(a+1)}\bv_{a(a+1)} + 2 c^{(b-1)b}\bv_{(b-1)b}  +2c^{b(b+1)}\bv_{b(b+1)}$
}
\begin{align*}
  & = (c^{aa}+c^{(a-1)a}+c^{a(a+1)})\bv_a + (c^{bb}+c^{(b-1)b}+c^{b(b+1)})\bv_b \\
  & \quad + c^{(a-1)a}\bv_{a-1}+ c^{a(a+1)}\bv_{a+1} + c^{(b-1)b}\bv_{b-1}+ c^{b(b+1)}\bv_{b+1}\\
  & = (c^{aa}+c^{(a-1)a}+c^{a(a+1)})\bv_a+ (c^{bb}+c^{(b-1)b}+c^{b(b+1)})\bv_b + sd_a\bv_a+ sd_b\bv_b \\
  & = (c^{aa}+s +sd_a) \bv_a+ (c^{bb}+s+sd_b)\bv_b.
\end{align*}
Substituting the fixed coordinates of $\bv_a = (-\ell,0)$ and $\bv_b=(\ell,0)$ reduces this expression to the vector
\[\left[\begin{array}{c}
(-c^{aa}-s-sd_a+c^{bb}+s+sd_b)\ell \\
0
\end{array}\right].\]

Finally, we address $\Qc$3.  Factoring the left side gives,
\flushleft{
$c^{aa}\bv_a\bv_a^T + c^{bb}\bv_b\bv_b^T + c^{(a-1)a}(\bv_{a-1}\bv_a^T+\bv_a\bv_{a-1}^T) + \cdots  + c^{b(b+1)}(\bv_b\bv_{b+1}^T+\bv_{b+1}\bv_b^T)$
}
\begin{align*}
 & = c^{aa}\bv_a\bv_a^T + c^{bb}\bv_b\bv_b^T \\
 & \qquad +(c^{(a-1)a}\bv_{a-1}+c^{a(a+1)}\bv_{a+1})\bv_a^T+ \bv_a(c^{(a-1)a}\bv_{a-1}^T+c^{a(a+1)}\bv_{a+1}^T) \\
 & \qquad +(c^{(b-1)b}\bv_{b-1}+c^{b(b+1)}\bv_{b+1})\bv_b^T+ \bv_b(c^{(b-1)b}\bv_{b-1}^T+c^{b(b+1)}\bv_{b+1}^T) \\
 & = c^{aa}\bv_a\bv_a^T + c^{bb}\bv_b\bv_b^T + sd_a\bv_a\bv_a^T + sd_b\bv_b\bv_b^T + \bv_a(sd_a\bv_a^T) + \bv_b(sd_b\bv_b^T) \\
 & = (c^{aa}+2sd_a)\bv_a\bv_a^T + (c^{bb}+2sd_b)\bv_b\bv_b^T.
\end{align*}
Again substituting the coordinates of $\bv_a$ and $\bv_b$, we obtain the matrix
\[\begin{bmatrix}
\left(c^{aa}+2sd_a +c^{bb}+2sd_b\right)\ell^2 & 0 \\
0 & 0
\end{bmatrix}.\]
The right side of $\Qc$3 is
\[\bv_a\bv_b^T+\bv_b\bv_a^T = \begin{bmatrix}
-2\ell^2 & 0 \\
0 & 0
\end{bmatrix}.\]
Hence the only equation that must be satisfied is exactly (\ref{eq:gen9}).
\end{proof}

\begin{remark}
We note that $s$ was specifically chosen so that (\ref{eq:gen7})-(\ref{eq:gen9}) would hold. The case $s=1$ happens when $d_a=-d_b$, i.e.\ only for the quadrilateral.
\end{remark}
\begin{theorem}
\label{thm:genpoly}
Given a convex polygon satisfying G\ref{g:ratio}, G\ref{g:minedge} and G\ref{g:maxangle}, $||\A||$ is uniformly bounded.
\end{theorem}
\begin{proof}
By Lemma~\ref{lem:cbd-implies-Abd}, it suffices to show a uniform bound on the six coefficients defined by equations (\ref{eq:gen1})-(\ref{eq:gen6}).  First we prove a uniform bound on $d_a$ and $d_b$ given G\ref{g:ratio}-G\ref{g:maxangle}.

\begin{figure}
\centering
\includegraphics[width=.4\textwidth]{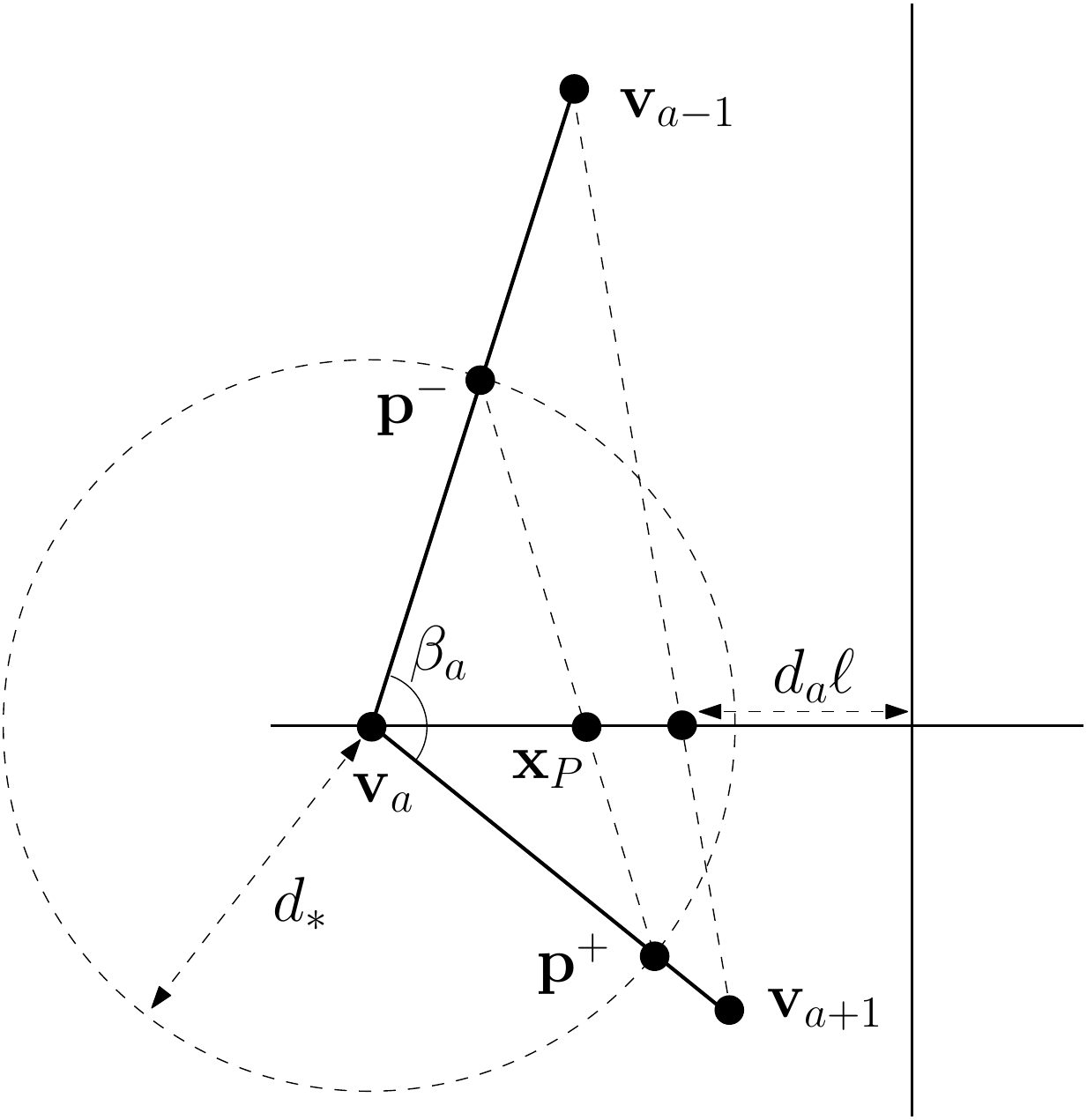}
\;\;
\includegraphics[width=.3\textwidth]{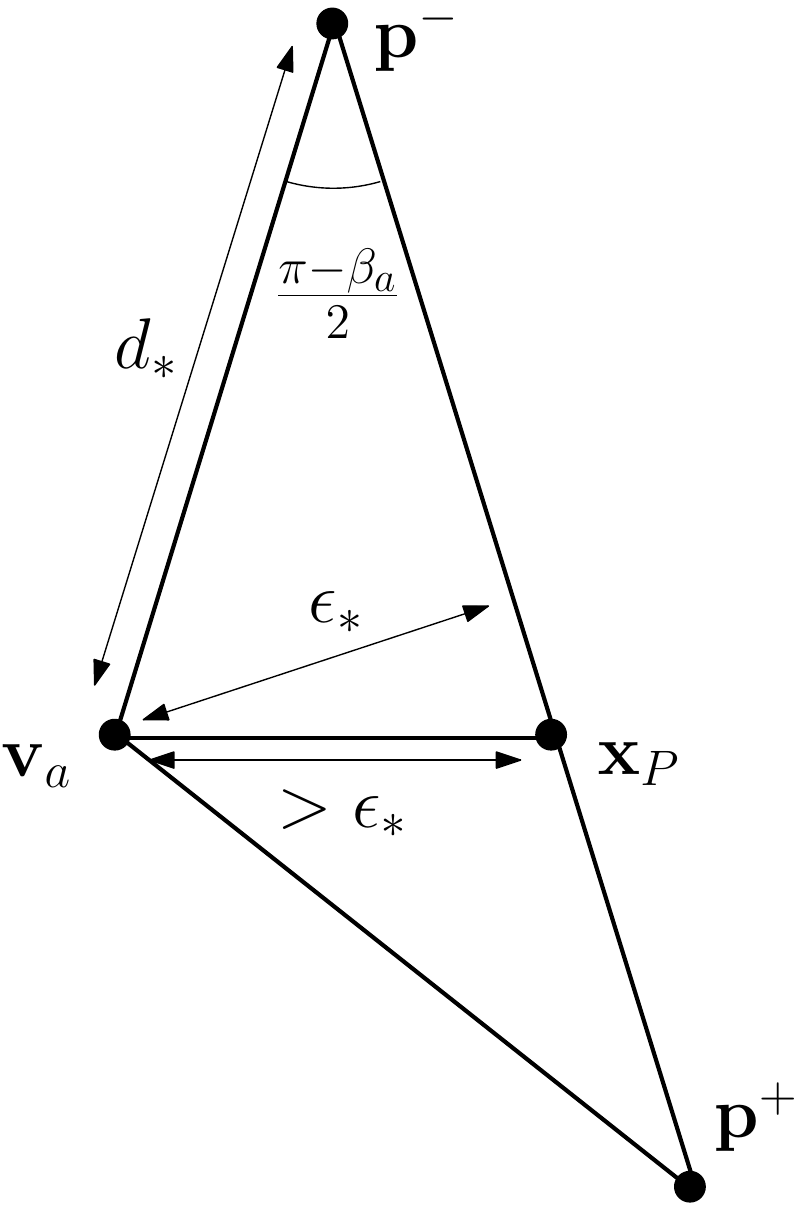}
\caption{Notation used in proof of Theorem~\ref{thm:genpoly}.}
\label{fg:proofGenPoly}
\end{figure}

We fix some notation as shown in Figure~\ref{fg:proofGenPoly}.  Let $C(\bv_a,d_\ast)$ be the circle of radius $d_\ast$ (from G\ref{g:minedge}) around $\bv_a$.  Let $\bp^-:=(p^-_x,p^-_y)$ and $\bp^+:=(p^+_x,p^+_y)$ be the points on $C(\bv_a,d_\ast)$ where the line segments to $\bv_a$ from $\bv_{a-1}$ and $\bv_{a+1}$, respectively, intersect.  The chord on $C(\bv_a,d_\ast)$ between $\bp^-$ and $\bp^+$ intersects the $x$-axis at $\bx_p:=(x_p,0)$.  By convexity, $(\bv_a)_x<x_p$.

To bound $x_p-(\bv_a)_x$ below, note that the triangle $\bv_a\bp^-\bp^+$ with angle $\beta_a$ at $\bv_a$ is isosceles.  Thus, the triangle $\bv_a\bp^-\bx_p$ has angle $\angle\bv_a\bp^-\bx_p=(\pi-\beta_a)/2$, as shown at the right of Figure~\ref{fg:proofGenPoly}.  
The distance to the nearest point on the line segment between $\bp^-$ and $\bp^+$ is $d_* \sin \left(\frac{\pi-\beta_a}{2} \right)$.
Based on G\ref{g:maxangle}, $\eps_\ast>0$ is defined to be
\begin{equation}
\label{eq:defepsast}
x_p-(\bv_a)_x \geq d_\ast\sin\left(\frac{\pi-\beta_a}{2}\right)> d_\ast\sin\left(\frac{\pi-\beta^\ast}{2}\right)=:\eps_\ast>0.
\end{equation}
Since $-d_a\ell < 1$ is the $x$-intercept of the line between $\bv_{a-1}$ and $\bv_{a+1}$, we have $x_p\leq -d_a\ell$.  Then we rewrite $(\bv_a)_x=-\ell$ in the geometrically suggestive form
\[(x_p-(\bv_a)_x)+(-d_a\ell-x_p)+(0+d_a\ell)=\ell.\]
Since $-d_a\ell-x_p\geq 0$, we have $\bx_p-(\bv_a)_x+d_a\ell\leq \ell$.  Using (\ref{eq:defepsast}), this becomes $d_a\ell<\ell-\eps_\ast$.  Recall from Figure~\ref{fg:genpoly} and previous discussion that $d_a,d_b\in[-1,1]$ and $-d_a\leq d_b$.  By symmetry, $d_b\ell<\ell-\eps_\ast$ and hence $d_a+d_b<2\ell-2\eps_\ast < 1-2\eps_\ast$.

We use the definition of $c^{aa}$ from (\ref{eq:gen5}), the derived bounds on $d_a$ and $d_b$, and the fact that $\ell\leq 1/2$ to conclude that
\[
|c^{aa}|<\frac{|2+2d_a|}{1+2\eps_\ast}<\frac{2+2(1-(\eps_\ast/\ell))}{1+2\eps_\ast}\leq \frac{4-4\eps_\ast}{1+2\eps_\ast}<4.
\]
Similarly, $|c^{bb}|<\frac{4-4\eps_\ast}{1+2\eps_\ast}<4$.  For the remaining coefficients, observe that the definition of $s$ in (\ref{eq:defofs}) implies that $0<s<2/(1+2\eps_\ast)$.  Equation (\ref{eq:gen1}) and the $y$-component of equation (\ref{eq:gen2}) ensure that $c^{(a-1)a}/s$ and $c^{a(a+1)}/s$ are the coefficients of a convex combination of $\bv_{a-1}$ and $\bv_{a+1}$.  Thus $c^{(a-1)a},c^{(a+1)a}\in (0,s)$ and $s$ serves as an upper bound on the norms of each coefficient.  Likewise, $|c^{(b-1)b}|,|c^{b(b+1)}|<s$.  Therefore,
\[
\max\left( \frac{4-4\eps_\ast}{1+2\eps_\ast}, \frac{2}{1+2\eps_\ast},1\right)
\]
is a uniform bound on all the coefficients of $\A$.
\end{proof}

\section{Converting Serendipity Elements to Lagrange-like Elements}
\label{sec:lag-reduc}

The $2n$ basis functions constructed thus far naturally correspond to vertices and edges of the polygon, but the functions associated to midpoints are not Lagrange-like.
This is due to the fact that functions of the form $\xi_{i(i+1)}$ may not evaluate to 1 at $\bv_{i(i+1)}$ or $\xi_{ii}$ may not evaluate to 0 at $\bv_{i(i+1)}$, even though the set of $\{\xi_{ij}\}$ satisfies the partition of unity property $\Qx$1.  
To fix this, we apply a simple bounded linear transformation given by the matrix $\B$ defined below.

\begin{figure}
\centering
\includegraphics[width=0.3\textwidth]{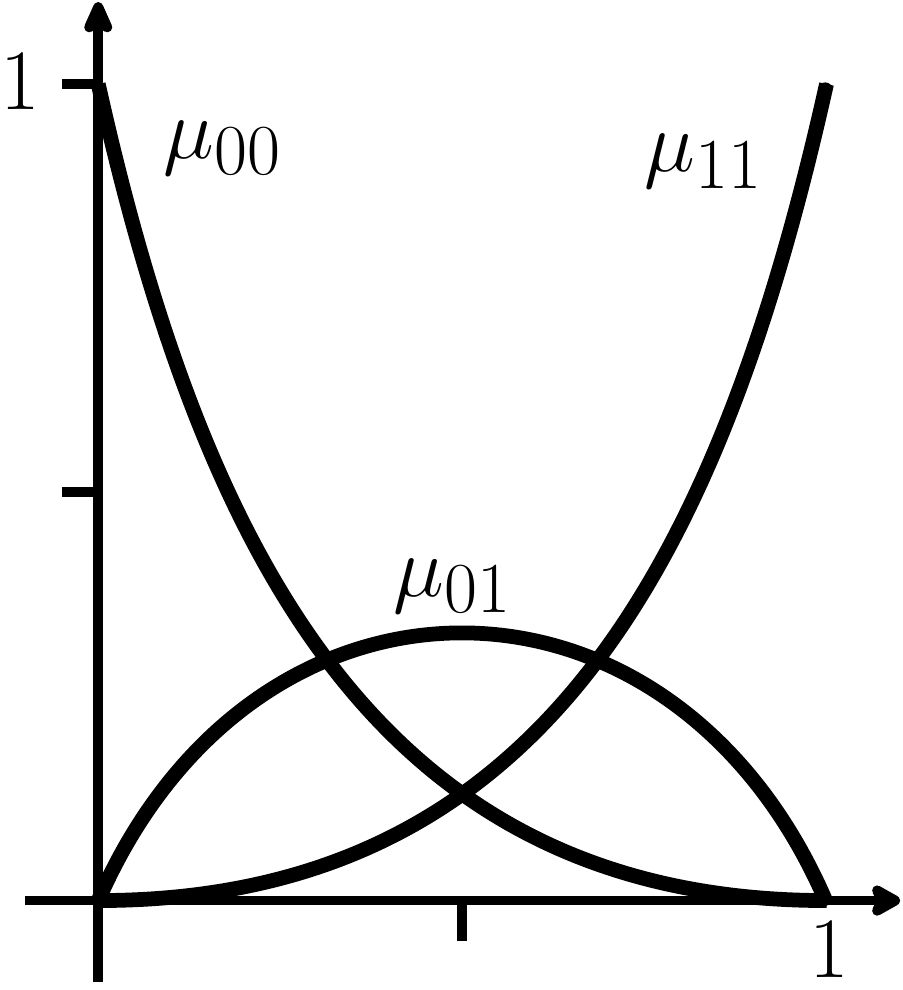}
\hspace{.2in}
\includegraphics[width=0.3\textwidth]{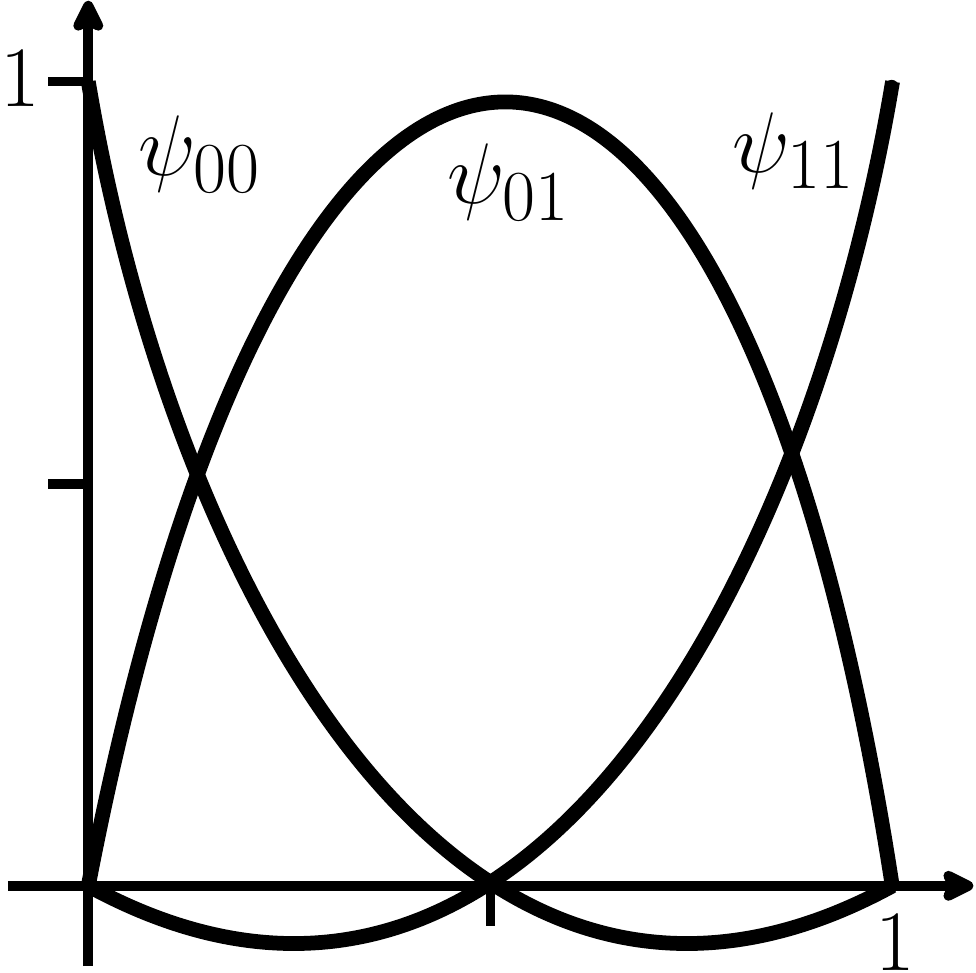}
\caption{A comparison of the product barycentric basis (left) with the standard Lagrange basis (right) for quadratic polynomials in one dimension.}
\label{fg:bases1d}
\end{figure}

To motivate our approach, we first consider a simpler setting: polynomial bases over the unit segment $[0,1]\subset\R$.  The barycentric functions on this domain are $\lambda_0(x) = 1 - x$, and $\lambda_1(x) =  x$.  Taking pairwise products, we get the quadratic basis $\mu_{00}(x) := \left(\lambda_0(x)\right)^2  = (1 - x)^2$, $\mu_{01}(x) := \lambda_0(x)\lambda_1(x) = (1 - x)x$, and $\mu_{11}(x) := \left(\lambda_1(x)\right)^2 = x^2$, shown on the left of Figure~\ref{fg:bases1d}.  This basis is not Lagrange-like since $\mu_{01}(1/2)\not=1$ and $\mu_{00}(1/2),\mu_{11}(1/2)\not=0$.  The quadratic Lagrange basis is given by $\psi_{00}(x) := 2(1-x)\left(\frac{1}{2}-x\right)$, $\psi_{01}(x) := 4(1-x)x$, and $\psi_{11}(x) := 2\left(x - \frac{1}{2}\right)x$, shown on the right of Figure~\ref{fg:bases1d}.  These two bases are related by the linear transformation $\B_{\text{1D}}$:
\begin{align}
[\psi_{ij}] =
\begin{bmatrix}\psi_{00} \\ \psi_{11} \\ \psi_{01} \end{bmatrix}
 &=
\begin{bmatrix}
1 & 0 & -1 \\
0 & 1 & -1 \\
0 & 0 & 4
\end{bmatrix}
\begin{bmatrix} \mu_{00}\\ \mu_{11} \\ \mu_{01}\end{bmatrix}
= \B_{\text{1D}} [\mu_{ij}].
\label{eq:barytolagrange}
\end{align}

This procedure generalizes to the case of converting the 2D serendipity basis $\{\xi_{ij}\}$ to a Lagrange like basis $\{\psi_{ij}\}$.  
Define
\[\psi_{ii} := \xi_{ii} - \xi_{i,i+1} - \xi_{i-1,i}\quad\text{and}\quad\psi_{i,i+1} = 4~\xi_{i,i+1}.\]
Using our conventions for basis ordering and index notation, the transformation matrix $\B$ taking $[\xi_{ij}]$ to $[\psi_{ij}]$ has the structure
\[
[\psi_{ij}] =
\begin{bmatrix}
\psi_{11} \\ \psi_{22} \\ \vdots \\ \psi_{nn} \\
\psi_{12} \\ \psi_{23} \\ \vdots \\ \psi_{n1}
\end{bmatrix}
 =
\left [
\begin{array}{ccccc|ccccc}
1 &   &   &   &   & -1 &    & \cdots &  & -1 \\
  & 1 &   &   &   & -1 & -1 & \cdots &    \\
  &   &\ddots & & &    & \ddots & \ddots\\
  &   & &\ddots & &    &    & \ddots & \ddots\\
  &   &   &   & 1 &    &    &   & -1 & -1 \\
\hline
  &&   &   &   &  4 \\
  &&   &   &   &  & 4 \\
  && 0 &   &   &  && \ddots \\
  &&   &   &   &  &&& \ddots \\
  &&   &   &   &  &&&& 4
\end{array}
\right]
\begin{bmatrix}
\xi_{11} \\ \xi_{22} \\ \vdots \\ \xi_{nn} \\
\xi_{12} \\ \xi_{23} \\ \vdots \\ \xi_{n1}
\end{bmatrix}
= \B [\xi_{ij}].
\]
The following proposition says that the functions $\{\psi_{ij}\}$ defined by the above transformation are Lagrange-like.
\begin{proposition}
For all $i,j \in \{1,\ldots,n\}$,
$\psi_{ii}(\bv_j) = \delta_i^j$,
$\psi_{ii}(\bv_{j,j+1}) = 0$,
$\psi_{i(i+1)}(\bv_j) = 0$, and
$\psi_{i(i+1)}(\bv_{j,j+1}) = \delta_i^j$.
\end{proposition}
\begin{proof}
We show the last claim first.  By the definitions of $\B$ and $\A$, we have
\begin{align*}
\psi_{i(i+1)}(\bv_{j,j+1}) 
& = 4~\xi_{i(i+1)}(\bv_{j,j+1}) = 4\left(\sum_{a=1}^n c_{aa}^{i(i+1)}\mu_{aa}(\bv_{j,j+1}) + \sum_{a<b} c_{ab}^{i(i+1)}\mu_{ab}(\bv_{j,j+1})\right)
\end{align*}
Since $\lambda_j$ is piecewise linear on the boundary of the polygon, $\lambda_a(\bv_{j,j+1})=1/2$ if $a\in\{j,j+1\}$ and zero otherwise.
Accordingly, $\mu_{aa}(\bv_{j,j+1})= 1/4$ if $a\in\{j,j+1\}$ and zero otherwise, while $\mu_{ab}(\bv_{j,j+1})=1/4$ if $\{a,b\}=\{j,j+1\}$ and zero otherwise.
\begin{align*}
\psi_{i(i+1)}(\bv_{j,j+1}) 
& = 4\left(\left(c_{jj}^{i(i+1)}+c_{(j+1)(j+1)}^{i(i+1)}\right)\cdot\frac 14 + c_{j(j+1)}^{i(i+1)}\cdot \frac 14\right) = c_{j(j+1)}^{i(i+1)} =  \delta_{ij},
\end{align*}
since the identity structure of $\A$ as given in (\ref{eq:A-struc}) implies that $c_{jj}^{i(i+1)}=c_{(j+1)(j+1)}^{i(i+1)}=0$ and that $c_{j(j+1)}^{i(i+1)} =  \delta_{ij}$.

Next, observe that $\mu_{ab}(\bv_j)=\lambda_a(\bv_j)\lambda_b(\bv_j)=1$ if $a=b=j$ and 0 otherwise.
Hence, any term of the form $c_{ab}^{**}\mu_{ab}(\bv_j)$ for $a\not=b$ is necessarily zero.
Therefore, by a similar expansion, $\psi_{i(i+1)}(\bv_j)=c_{jj}^{i(i+1)}=0$, proving the penultimate claim.

For the first two claims, similar analysis yields
\begin{align*}
\psi_{ii}(\bv_j) 
& = \xi_{ii}(\bv_j) - \xi_{i(i+1)}(\bv_j) - \xi_{(i-1)i}(\bv_j) \\
& = c_{jj}^{ii}\cdot 1 -  c_{jj}^{i(i+1)}\cdot 1 - c_{jj}^{(i-1)i}\cdot 1\\
& = c_{jj}^{ii} =  \delta_{ij},
\end{align*}
again by the identity structure of $\A$.
Finally, by similar analysis, we have that
\begin{align*}
\psi_{ii}(\bv_{j,j+1}) 
& = \xi_{ii}(\bv_{j,j+1}) - \xi_{i(i+1)}(\bv_{j,j+1}) - \xi_{(i-1)i}(\bv_{j,j+1}) \\
& = (c_{jj}^{ii} +  c_{(j+1)(j+1)}^{ii} + c_{j(j+1)}^{ii})\frac 14  - \xi_{i(i+1)}(\bv_{j,j+1}) - \xi_{(i-1)i}(\bv_{j,j+1}) \\
& = (\delta_{ij} +\delta_{i(j+1)})\frac 14  - \frac 14\delta_{ij} - \frac 14\delta_{i(j+1)} = 0,
\end{align*}
completing the proof.
\end{proof}
In closing, note that $||\B||$ is uniformly bounded since its entries all lie in $\{-1,0,1,4\}$.

\section{Applications and Extensions}
\label{sec:conc}

\begin{figure}[b]
\centering
\begin{tabular}{ccc}
\parbox{.24\textwidth}{\includegraphics[width=.23\textwidth]{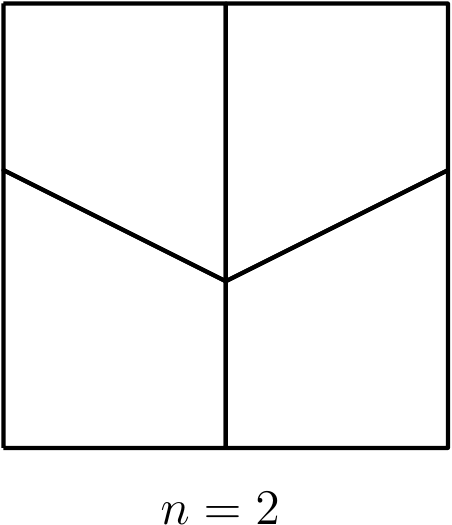}}
&
\parbox{.24\textwidth}{\includegraphics[width=.23\textwidth]{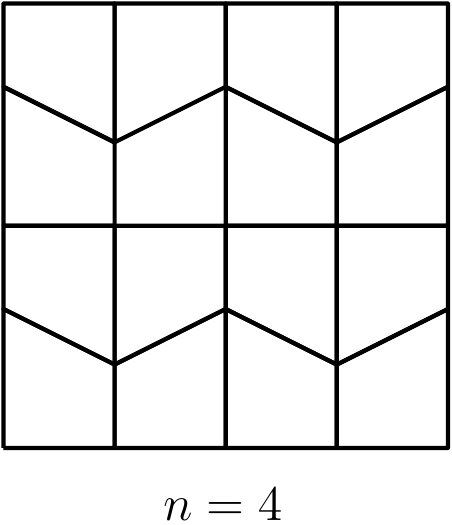}}
&
{\small 
\begin{tabular}{c|cc|cc}
 & \multicolumn{2}{c|}{$\vn{u-u_h}_{L^2}$} & \multicolumn{2}{c}{$\vn{\nabla(u-u_h)}_{L^2}$}  \\ \hline
n & error & rate & error & rate \\ \hline
2 & 2.34e-3 & & 2.22e-2 & \\
4 & 3.03e-4 & 2.95& 6.10e-3 & 1.87\\
8 & 3.87e-5 & 2.97& 1.59e-3 & 1.94\\
16 & 4.88e-6 & 2.99& 4.04e-4 & 1.97\\
32 & 6.13e-7& 3.00& 1.02e-4 & 1.99\\
64 & 7.67e-8& 3.00& 2.56e-5 & 1.99\\
128 & 9.59e-9& 3.00& 6.40e-6 & 2.00\\
256 & 1.20e-9& 3.00& 1.64e-6 & 1.96\\
\end{tabular} }
\end{tabular}
\caption{Trapezoidal meshes (left) fail to produce quadratic convergence with traditional serendipity elements; see~\cite{ABF02}.  
Since our construction begins with affinely-invariant generalized barycentric functions, the expected quadratic convergence rate can be recovered (right).  
The results shown were generated using the basis $\{\psi_{ij} \}$ resulting from the selection of the mean value coordinates as the initial barycentric functions.  
}
\label{fg:nonaffinequad}
\end{figure}

Our quadratic serendipity element construction has a number of uses in modern finite element application contexts.  
First, the construction for quadrilaterals given in Section~\ref{subsec:quads} allows for quadratic order methods on \textit{arbitrary} quadrilateral meshes with only eight basis functions per element instead of the nine used in a bilinear map of the biquadratic tensor product basis on a square.
In particular, we show that our approach maintains quadratic convergence on a mesh of convex quadrilaterals known to result in only linear convergence when traditional serendipity elements are mapped non-affinely~\cite{ABF02}.

We solve Poisson's equation on a square domain composed of  $n^2$ trapezoidal elements as shown in Figure~\ref{fg:nonaffinequad}~(left).  
Boundary conditions are prescribed according to the solution $u(x,y) = \sin(x) e^y$; we use our construction from Section~\ref{subsec:quads} starting with mean value coordinates $\{\lmval_i\}$. 
Mean value coordinates were selected based on a few advantages they have over other types: they are easy to compute based on an explicit formula and the coordinate gradients do not degrade based on large interior angles~\cite{RGB2011b}.
For this particular example, where no interior angles asymptotically approach $180^\circ$, Wachspress coordinates give very similar results. 
As shown in Figure~\ref{fg:nonaffinequad}~(right), the expected convergence rates from our theoretical analysis are observed, namely, cubic in the $L^2$-norm and quadratic in the $H^1$-norm.

\begin{figure}[ht]
\centering
\includegraphics[width=.3\textwidth]{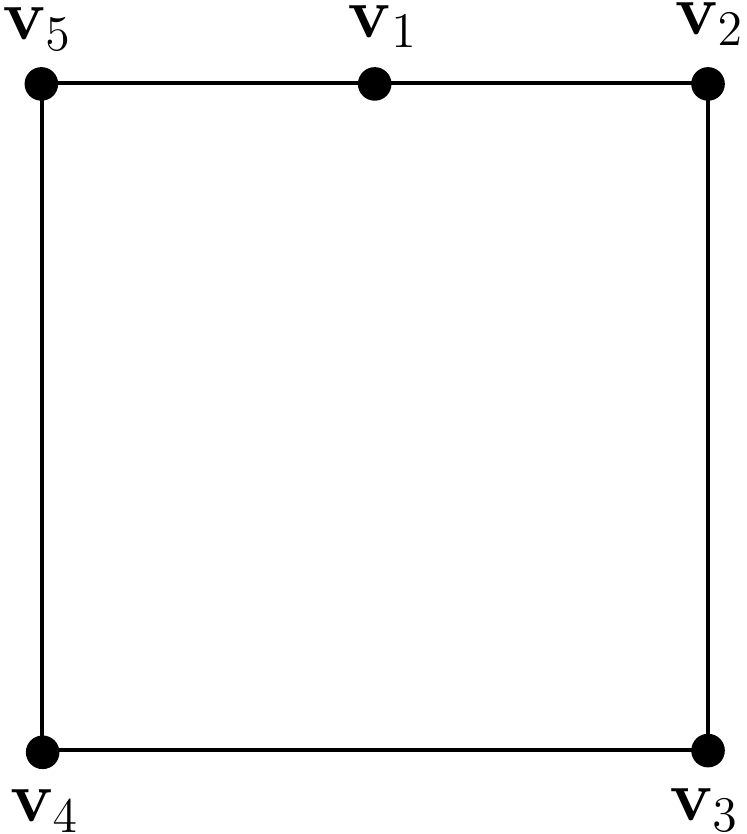}
\;\;
\includegraphics[width=.3\textwidth]{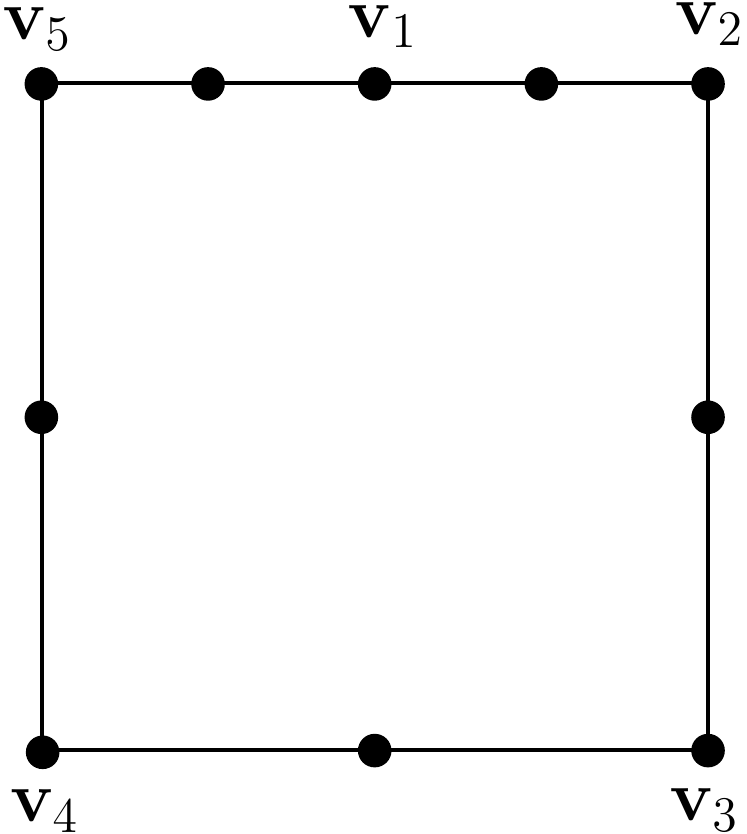}
\caption{
Theorem~\ref{thm:genpoly} can be generalized to allow certain types of geometries that do not satisfy G\ref{g:maxangle}.
The degenerate pentagon (left), widely used in adaptive finite element methods for quadrilateral meshes, satisfies G\ref{g:ratio} and  G\ref{g:minedge}, but only satisfies G\ref{g:maxangle} for four of its vertices.  
The bounds on the coefficients $c^{ij}_{ab}$ from Section~\ref{sec:genpoly} still hold on this geometry, resulting in the Lagrange-like quadratic element (right).}
\label{fg:degenPent}
\end{figure}

An additional application of our method is to adaptive finite elements, such as the one shown in Figure~\ref{fg:degenPent}.  
This is possible since the result of Theorem~\ref{thm:genpoly} still holds if G\ref{g:maxangle} fails to hold only on a set of consecutive vertices of the polygon.  
This weakened condition suffices since consecutive large angles in the polygon do not cause the coefficients $c^{ij}_{ab}$ to blow up.  
For instance, consider the degenerate pentagon shown in Figure~\ref{fg:degenPent} which satisfies this weaker condition but not G\ref{g:maxangle}.  
Examining the potentially problematic coefficients $c^{ij}_{25}$, observe that the lines through $\bv_1$, $\bv_4$ and $\bv_1$, $\bv_3$ both intersect the midpoint of the line through $\bv_2$, $\bv_5$ (which happens to be $\bv_1$).
In the computation of the $c^{ij}_{25}$ coefficients, the associated values $d_2$ and $d_5$ are both zero and hence $s=1$~(recall Figure~\ref{fg:genpoly} and formula (\ref{eq:defofs})).
Since $s$ is bounded away from $\infty$, the analysis from the proof of Theorem~\ref{thm:genpoly} holds as stated for these coefficients and hence for the entire element. 
A more detailed analysis of such large-angle elements is an open question for future study.

\begin{figure}[ht]
\centering
\includegraphics[width=.4\textwidth]{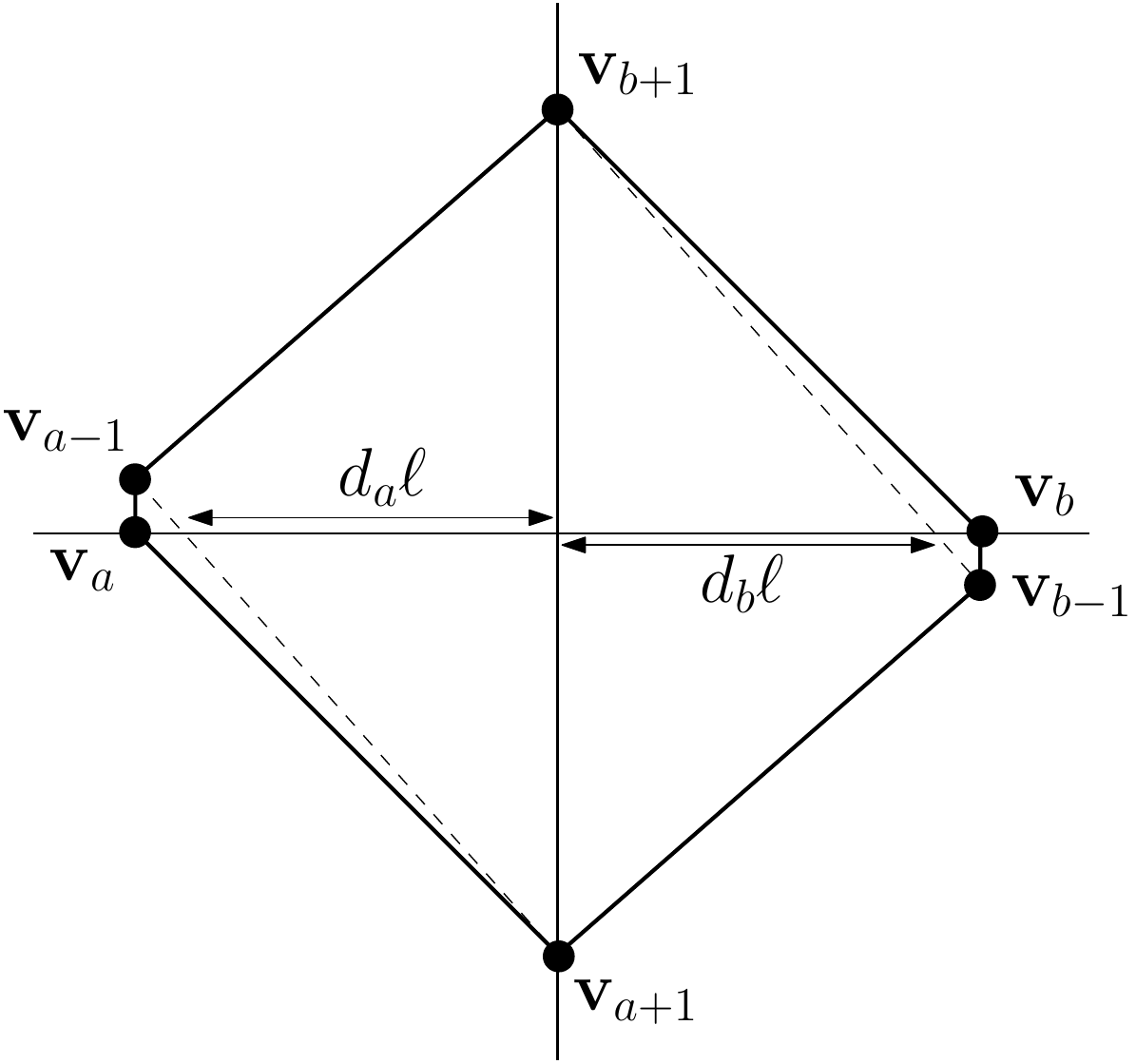}
\;\;
\includegraphics[width=.4\textwidth]{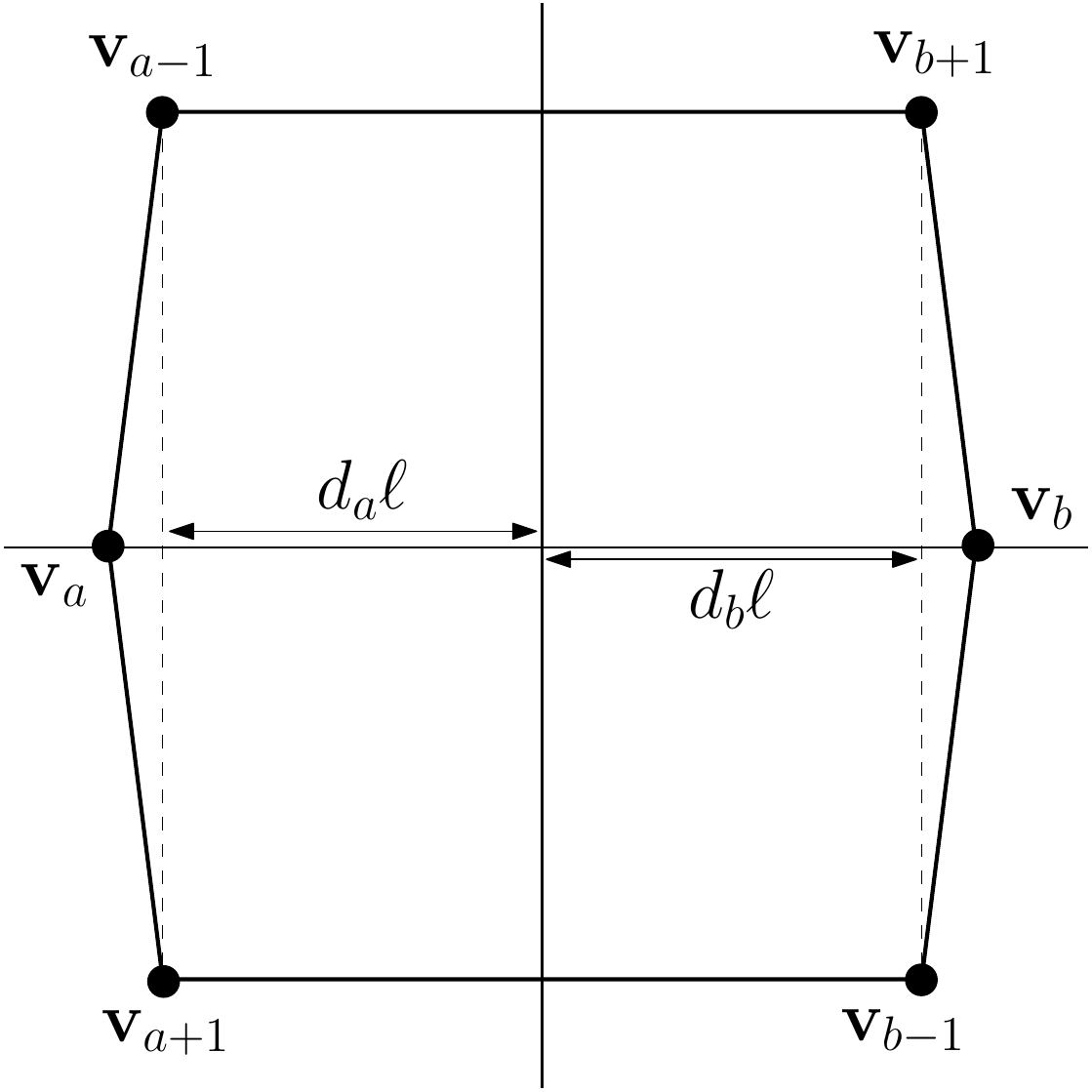}
\caption{
The hypotheses of Theorem~\ref{thm:genpoly} cannot be relaxed entirely as demonstrated by these shapes.
If G\ref{g:minedge} does not hold, arbitrarily small edges can cause a blowup  in the coefficients $c^{ij}_{ab}$ (left).
If G\ref{g:maxangle} does not hold, non-consecutive angles approaching $\pi$ can cause a similar blowup.}
\label{fg:geomNecc}
\end{figure}

Nevertheless, the geometric hypotheses of Theorem~\ref{thm:genpoly} cannot be relaxed entirely.
Arbitrarily large non-consecutive large angles as well as very short edges,  can cause a blowup in the coefficients used in the construction of $\A$, as shown in Figure~\ref{fg:geomNecc}.
In the left figure, as edges $\bv_{a-1}\bv_a$ and $\bv_{b-1}\bv_b$ approach length zero, $d_a$ and $d_b$ both approach one, meaning $s$ (in the construction of Section~\ref{sec:genpoly}) approaches $\infty$.  
In this case, the coefficients $c^{(a-1)a}_{ij}$ and $c^{(b-1)b}_{ij}$ grow larger without bound, thereby violating the result of Theorem~\ref{thm:genpoly}.  
In the right figure, as the overall shape approaches a square, $d_a$ and $d_b$ again approach one so that $s$ again approaches $\infty$.  In this case, all the coefficients $c^{(a-1)a}_{ij}$, $c^{a(a+1)}_{ij}$, $c^{(b-1)b}_{ij}$ and $c^{b(b+1)}_{ij}$ all grow without bound.  
Nevertheless, if these types of extreme geometries are required, it may be possible to devise alternative definitions of the $c^{ab}_{ij}$ coefficients satisfying $\Qc$1-$\Qc$3 with controlled norm estimates since the set of restrictions $\Qc$1-$\Qc$3 does not have full rank. 
Note that this flexibility has lead to multiple constructions of the traditional serendipity square~\cite{MH92,KOF99}. 
Cursory numerical experimentation suggests that some bounded construction exists even in the degenerate situation.  

The computational cost of our method is an important consideration to application contexts.
A typical finite element method using our approach would involve the following steps: (1) selecting $\lambda_i$ coordinates and implementing the corresponding $\psi_{ij}$ basis functions, (2) defining a quadrature rule for each affine-equivalent class of shapes appearing in the domain mesh, (3) assembling a matrix $\Ell$ representing the discrete version of linear operator, and (4) solving a linear system of the form $\Ell u = f$.
The quadrature step may incur some computational effort, however, if only a few shape templates are needed, this is a one-time fixed pre-preprocessing cost.
In the trapezoidal mesh example from Figure~\ref{fg:nonaffinequad}, for instance, we only needed one quadrature rule as all domain shapes were affinity equivalent.
Assembling the matrix $\Ell$ may also be expensive as the entries involve integrals of products of gradients of $\psi_{ij}$ functions.
Again, however, this cost is incurred only once per affine-equivalent domain shape and thus can be reasonable to allow, depending on the application context.

The computational advantage to our approach comes in the final linear solve.  
The size of the matrix $\Ell$ is proportional to the number of edges in the mesh, matching the size of the corresponding matrix for quadratic Lagrange elements on triangles or quadratic serendipity elements on squares.
If the pairwise products $\mu_{ab}$ were used instead of the $\psi_{ij}$ functions, the size of $\Ell$ would be proportional to the \textit{square} of the number of edges in the mesh, a substantial difference.

Finally, we note that although this construction is specific to quadratic elements, the approach seems adaptable, with some effort, to the construction of cubic and higher order serendipity elements on generic convex polygons. 
As a larger linear system must be satisfied, stating an explicit solution becomes complex.  
Further research along these lines should probably assert the existence of a uniformly bounded solution without specifying the construction. 
In practice, a least squares solver could be used to construct such a basis numerically.

\bibliographystyle{amsplain}
\bibliography{../references}

\providecommand{\bysame}{\leavevmode\hbox to3em{\hrulefill}\thinspace}
\providecommand{\MR}{\relax\ifhmode\unskip\space\fi MR }
\providecommand{\MRhref}[2]{%
  \href{http://www.ams.org/mathscinet-getitem?mr=#1}{#2}
}
\providecommand{\href}[2]{#2}
\begin{thebibliography}{10}

\bibitem{apprato1979rational}
D.~Apprato, R.~Arcanceli, and J.~L. Gout, \emph{{Rational interpolation of
  Wachspress error estimates}}, Comput. Math. Appl. \textbf{5} (1979), no.~4,
  329--336.

\bibitem{apprato1979elements}
D.~Apprato, R.~Arcangeli, and J.~L. Gout, \emph{Sur les elements finis
  rationnels de {W}achspress}, Numer. Math. \textbf{32} (1979), no.~3,
  247--270.

\bibitem{AA10}
D.~N. Arnold and G.~Awanou, \emph{The serendipity family of finite elements},
  Found. Comput. Math. \textbf{11} (2011), no.~3, 337--344.

\bibitem{ABF02}
D.~N. Arnold, D.~Boffi, and R.~S. Falk, \emph{{Approximation by quadrilateral
  finite elements}}, Math. Comput. \textbf{71} (2002), no.~239, 909--922.

\bibitem{BS08}
S.~C. Brenner and L.~R. Scott, \emph{The mathematical theory of finite element
  methods}, third ed., Texts in Applied Mathematics, vol.~15, Springer, New
  York, 2008. \MR{2373954 (2008m:65001)}

\bibitem{C1989}
D.~Chavey, \emph{Tilings by regular polygons {II}: A catalog of tilings},
  Comput. Math. Appl. \textbf{17} (1989), no.~1--3, 147--165.

\bibitem{C2008}
S.~H. Christiansen, \emph{A construction of spaces of compatible differential
  forms on cellular complexes}, Math. Models Methods Appl. Sci. \textbf{18}
  (2008), no.~5, 739--757.

\bibitem{Ci02}
P.~G. Ciarlet, \emph{The finite element method for elliptic problems}, second
  ed., Classics in Applied Mathematics, vol.~40, SIAM, Philadelphia, PA, 2002.

\bibitem{CSCMCD2003}
E.~Cueto, N.~Sukumar, B.~Calvo, M.~A. Mart{\'{\i}}nez, J.~Cego{\~n}ino, and
  M.~Doblar{\'e}, \emph{Overview and recent advances in natural neighbour
  {G}alerkin methods}, Arch. Comput. Methods Engrg. \textbf{10} (2003), no.~4,
  307--384. \MR{2032470 (2004m:65190)}

\bibitem{DL04}
S.~Dekel and D.~Leviatan, \emph{{The Bramble-Hilbert lemma for convex
  domains}}, SIAM J. Math. Anal. \textbf{35} (2004), no.~5, 1203--1212.

\bibitem{EG04}
A.~Ern and J.-L. Guermond, \emph{Theory and practice of finite elements},
  Applied Mathematical Sciences, vol. 159, Springer-Verlag, New York, 2004.
  \MR{2050138 (2005d:65002)}

\bibitem{F1990}
G.~Farin, \emph{{Surfaces over Dirichlet tessellations}}, Comput. Aided Geom.
  Des. \textbf{7} (1990), no.~1-4, 281--292.

\bibitem{F2003}
M.~Floater, \emph{Mean value coordinates}, Comput. Aided Geom. Des. \textbf{20}
  (2003), no.~1, 19--27.

\bibitem{FHK2006}
M.~Floater, K.~Hormann, and G.~K{\'o}s, \emph{A general construction of
  barycentric coordinates over convex polygons}, Adv. Comput. Math. \textbf{24}
  (2006), no.~1, 311--331.

\bibitem{GB2010}
A.~Gillette and C.~Bajaj, \emph{A generalization for stable mixed finite
  elements}, Proc. 14th ACM Symp. Solid Phys. Modeling, 2010, pp.~41--50.

\bibitem{GB2011}
\bysame, \emph{Dual formulations of mixed finite element methods with
  applications}, Comput. Aided Des. \textbf{43} (2011), no.~10, 1213--1221,
  arXiv:1012.3929 [math.DG].

\bibitem{GRB2011}
A.~Gillette, A.~Rand, and C.~Bajaj, \emph{Error estimates for generalized
  barycentric coordinates}, to appear in Advances in Computational Mathematics
  (2011), 1--23.

\bibitem{gout1979construction}
J.~L. Gout, \emph{{Construction of a Hermite rational ``Wachspress type''
  finite element}}, Comput. Math. Appl. \textbf{5} (1979), no.~4, 337--347.

\bibitem{gout1985rational}
\bysame, \emph{{Rational Wachspress-type finite elements on regular hexagons}},
  IMA J. Numer. Anal. \textbf{5} (1985), no.~1, 59.

\bibitem{H1987}
T.~Hughes, \emph{The finite element method}, Prentice Hall Inc., Englewood
  Cliffs, NJ, 1987, Linear static and dynamic finite element analysis, With the
  collaboration of Robert M. Ferencz and Arthur M. Raefsky.

\bibitem{JMRGS07}
P.~Joshi, M.~Meyer, T.~DeRose, B.~Green, and T.~Sanocki, \emph{{Harmonic
  coordinates for character articulation}}, ACM Trans. Graphics \textbf{26}
  (2007), 71.

\bibitem{KOF99}
F.~Kikuchi, M.~Okabe, and H.~Fujio, \emph{Modification of the 8-node
  serendipity element}, Comput. Methods Appl. Mech. Engrg. \textbf{179} (1999),
  no.~1-2, 91--109.

\bibitem{LS08}
T.~Langer and H.P. Seidel, \emph{Higher order barycentric coordinates}, Comput.
  Graphics Forum, vol.~27, Wiley Online Library, 2008, pp.~459--466.

\bibitem{MH92}
R.~H. Mac{N}eal and R.~L. Harder, \emph{Eight nodes or nine?}, Int. J. Numer.
  Methods Eng. \textbf{33} (1992), no.~5, 1049--1058.

\bibitem{MKBWG2008}
S.~Martin, P.~Kaufmann, M.~Botsch, M.~Wicke, and M.~Gross, \emph{Polyhedral
  finite elements using harmonic basis functions}, Proc. Symp. Geom. Proc.,
  2008, pp.~1521--1529.

\bibitem{MP2008}
P.~Milbradt and T.~Pick, \emph{Polytope finite elements}, Int. J. Numer.
  Methods Eng. \textbf{73} (2008), no.~12, 1811--1835.

\bibitem{RGB2011b}
A.~Rand, A.~Gillette, and C.~Bajaj, \emph{Interpolation error estimates for
  mean value coordinates}, arXiv:1111.5588 (2011), 1--18.

\bibitem{RS2006}
M.~M. Rashid and M.~Selimotic, \emph{{A three-dimensional finite element method
  with arbitrary polyhedral elements}}, Int. J. Numer. Methods Eng. \textbf{67}
  (2006), no.~2, 226--252.

\bibitem{S1980}
R.~Sibson, \emph{A vector identity for the {D}irichlet tessellation}, Math.
  Proc. Cambridge Philos. Soc. \textbf{87} (1980), no.~1, 151--155.

\bibitem{SAB2010}
D.~Sieger, P.~Alliez, and M.~Botsch, \emph{Optimizing voronoi diagrams for
  polygonal finite element computations}, Proc. 19th Int. Meshing Roundtable
  (2010), 335--350.

\bibitem{SF73}
G.~Strang and G.~J. Fix, \emph{An analysis of the finite element method},
  Prentice-Hall Inc., Englewood Cliffs, N. J., 1973, Prentice-Hall Series in
  Automatic Computation.

\bibitem{SM2006}
N.~Sukumar and E.~A. Malsch, \emph{Recent advances in the construction of
  polygonal finite element interpolants}, Archives Comput. Methods. Eng.
  \textbf{13} (2006), no.~1, 129--163.

\bibitem{ST2004}
N.~Sukumar and A.~Tabarraei, \emph{{Conforming polygonal finite elements}},
  Int. J. Numer. Methods Eng. \textbf{61} (2004), no.~12, 2045--2066.

\bibitem{TS06}
A.~Tabarraei and N.~Sukumar, \emph{{Application of polygonal finite elements in
  linear elasticity}}, International Journal of Computational Methods
  \textbf{3} (2006), no.~4, 503--520.

\bibitem{Ve99}
R.~Verf{\"u}rth, \emph{A note on polynomial approximation in {S}obolev spaces},
  Math. Modelling Numer. Anal. \textbf{33} (1999), no.~4, 715--719.

\bibitem{W1975}
E.~L. Wachspress, \emph{A rational finite element basis}, Mathematics in
  Science and Engineering, vol. 114, Academic Press, New York, 1975.

\bibitem{WBG07}
M.~Wicke, M.~Botsch, and M.~Gross, \emph{{A finite element method on convex
  polyhedra}}, Comput. Graphics Forum \textbf{26} (2007), no.~3, 355--364.

\bibitem{ZK00}
J.~Zhang and F.~Kikuchi, \emph{Interpolation error estimates of a modified
  8-node serendipity finite element}, Numer. Math. \textbf{85} (2000), no.~3,
  503--524.

\bibitem{ZT2000}
O.~Zienkiewicz and R.~Taylor, \emph{The finite element method}, fifth ed.,
  Butterworth-Heinemann, London, 2000.

\end{thebibliography}

\end{document}